\documentclass{amsart}

\usepackage{amsmath,amsfonts,amssymb,enumitem}
\usepackage{tikz}

\newcommand{\thmcolor}{\color{black}}
\newcommand{\defcolor}{\color{black}}
\newcommand{\newcolor}{\color{blue!50!black}}

\newtheorem{theorem}{\thmcolor Theorem}[subsection]
\newtheorem{lemma}[theorem]{\thmcolor Lemma}
\newtheorem{proposition}[theorem]{\thmcolor Proposition}
\newtheorem{corollary}[theorem]{\thmcolor Corollary}
\newtheorem{conjecture}[theorem]{Conjecture}
\newtheorem*{theorem*}{Theorem}
\newtheorem*{conjecture*}{Conjecture}

\theoremstyle{definition}
\newtheorem{definition}[theorem]{\defcolor Definition}

\theoremstyle{remark}

\numberwithin{equation}{section}

\newcommand{\newword}[1]{{\newcolor\textbf{\emph{#1}}}}

\newcommand{\key}{\ensuremath{\kappa}}
\newcommand{\pinp}{\ensuremath{\mathcal{P}}}
\newcommand{\atom}{\ensuremath{\mathcal{A}}}

\newcommand{\SSRT}{\ensuremath{\mathrm{SSRT}}}

\newcommand{\KD}{\ensuremath{\mathrm{KD}}}
\newcommand{\PKD}{\ensuremath{\mathrm{PKD}}}
\newcommand{\AKD}{\ensuremath{\mathrm{AKD}}}

\newcommand{\LAT}{\ensuremath{\mathrm{LAT}}}
\newcommand{\LKT}{\ensuremath{\mathrm{LKT}}}

\newcommand{\D}{\ensuremath{\mathbb{D}}}

\newcommand{\wt}{\ensuremath\mathrm{\mathbf{wt}}}

\newcommand{\thread}{\ensuremath\boldsymbol{\theta}}
\newcommand{\e}{\ensuremath\mathrm{\mathbf{e}}}

\newcommand{\Rect}{\ensuremath\mathfrak{Rect}}
\newcommand{\Label}{\ensuremath\mathcal{L}}

\newcommand{\sort}{\ensuremath\boldsymbol{\lambda}}

\newcommand{\lswap}{\ensuremath\preceq}
\newcommand{\lswappin}{\ensuremath\lswap^{0}}

\newcommand{\kup}{\subseteq}
\newcommand{\kupdot}{\mathrel{\ooalign{$\subset$\cr\hidewidth\hbox{$\cdot$}\cr}}}

\newcommand{\kplus}{+^{\hspace{-0.5ex}(k)}}

\newcommand{\rkey}{\ensuremath\mathcal{R}}

\newlength\cellsize 

\savebox2{
\begin{picture}(8,8)
\put(0,0){\line(1,0){8}}
\put(0,0){\line(0,1){8}}
\put(8,0){\line(0,1){8}}
\put(0,8){\line(1,0){8}}
\end{picture}}

\newcommand\boxify[1]{\def\thearg{#1}\def\nothing{}%
\ifx\thearg\nothing\vrule width0pt height\cellsize depth0pt%
  \else\hbox to 0pt{\usebox2\hss}\fi%
  \vbox to \cellsize{\vss\hbox to \cellsize{\hss$_{#1}$\hss}\vss}}

\savebox3{
\begin{picture}(8,8)
\put(4,4){\circle{8}}
\end{picture}}

\newcommand{\circify}[1]{\def\thearg{#1}\def\nothing{}%
\ifx\thearg\nothing\vrule width0pt height\cellsize depth0pt%
  \else\hbox to 0pt{\usebox3\hss}\fi%
  \vbox to \cellsize{\vss\hbox to \cellsize{\hss$_{#1}$\hss}\vss}}

\newcommand\nullify[1]{\def\thearg{#1}\def\nothing{}%
\ifx\thearg\nothing\vrule width0pt height\cellsize depth0pt%
  \else\hbox to 0pt{\hss}\fi%
  \vbox to \cellsize{\vss\hbox to \cellsize{\hss$_{#1}$\hss}\vss}}

\newcommand\tableau[1]{\vtop{\let\\=\cr
\setlength\baselineskip{-8000pt}
\setlength\lineskiplimit{8000pt}
\setlength\lineskip{0pt}
\halign{&\boxify{##}\cr#1\crcr}}}

\newcommand\cirtab[1]{\vline\vtop{\let\\=\cr
\setlength\baselineskip{-8000pt}
\setlength\lineskiplimit{8000pt}
\setlength\lineskip{0pt}
\halign{&\circify{##}\cr#1\crcr}}}

\newcommand\nulltab[1]{\vtop{\let\\=\cr
\setlength\baselineskip{-8000pt}
\setlength\lineskiplimit{8000pt}
\setlength\lineskip{0pt}
\halign{&\nullify{##}\cr#1\crcr}}}

\savebox2{
\begin{picture}(8,8)
\put(0,0){\line(1,0){8}}
\put(0,0){\line(0,1){8}}
\put(8,0){\line(0,1){8}}
\put(0,8){\line(1,0){8}}
\end{picture}}

\newcommand{\csix}{red}
\newcommand{\cfiv}{orange}
\newcommand{\cfou}{yellow}
\newcommand{\cthr}{green}
\newcommand{\ctwo}{blue}
\newcommand{\cone}{violet}

\newcommand{\cball}[2]{%
  \begin{tikzpicture}
    \filldraw[fill=#1!35,draw=black] (0,0) circle (4pt) node {$\scriptstyle #2$};
  \end{tikzpicture}
}

\newcommand{\cbox}[2]{%
  \begin{tikzpicture}
    \draw +(-4pt,-4pt) rectangle +(4pt,4pt);
    \node (0,0) {$\scriptstyle \mathbf{#2}$};
  \end{tikzpicture}
}

\newcommand{\leftball}[2]{\makebox[0pt]{\raisebox{1.5pt}{$\leftarrow$}}\cball{#1}{#2}}

\usepackage{xcolor}
\usepackage[colorinlistoftodos]{todonotes}

\begin{document}


\title[Insertion algorithm for Demazure characters]{An insertion algorithm for multiplying Demazure characters by Schur polynomials}  

\author[S. Assaf]{Sami Assaf}
\address{Department of Mathematics, University of Southern California, 3620 S. Vermont Ave., Los Angeles, CA 90089-2532, U.S.A.}
\email{shassaf@usc.edu}






\begin{abstract}
  We introduce an insertion algorithm on Kohnert's combinatorial model for Demazure characters, generalizing Robinson--Schensted--Knuth insertion on tableaux. Our new insertion yields an explicit, nonnegative formula expressing the product of a Demazure character and a Schur polynomial as a sum of Schubert characters, partially resolving Polo's conjecture that the tensor product of Demazure modules admits a Schubert filtration. 
\end{abstract}

\maketitle

%
\section{Introduction}
%
\label{sec:introduction}

For $G$ a connected, simply-connected, semi-simple Lie group, the finite dimensional simple $G$-modules $V(\lambda)$ are indexed by dominant weights $\lambda$. For $w \in W$ the Weyl group, the \emph{Demazure module} $V_w(\lambda)$ is the $B$-submodule of $V(\lambda)$ generated by the one-dimensional extremal weight space of weight $w \lambda$. Bruhat order filters Demazure modules by $V_u(\lambda) \subseteq V_w(\lambda)$ if and only if $u \preceq w$. Thus Demazure modules interpolate between the one-dimensional highest weight space $V_{\mathrm{id}}(\lambda)$ and the irreducible module $V(\lambda) = V_{w_0}(\lambda)$, for $w_0$ the long element. Since $V_u(\lambda) = V_w(\lambda)$ if and only if $u \lambda = w \lambda$, we index Demazure modules by the extremal weight $w \lambda$.

For $G= GL_n(\mathbb{C})$ the general linear group, the weights are weak compositions $\alpha = w \lambda \in \mathbb{Z}^n_{\ge 0}$, and the \emph{Demazure characters} $\key_{\alpha} = \mathrm{ch}(V_w(\lambda))$ form a $\mathbb{Z}$-basis of the polynomial ring $\mathbb{Z}[x_1,\ldots,x_n]$. Thus it is natural to consider structure constants $c_{\beta,\alpha}^{\gamma}\in\mathbb{Z}$ for Demazure characters defined by
\begin{equation}
  \key_{\beta} \key_{\alpha} = \sum_{\gamma} c_{\beta,\alpha}^{\gamma} \key_{\gamma}.
  \label{e:LRC-dem}
\end{equation}

Mathieu \cite{Mat89} proved a conjecture of Polo \cite{Pol89} that if one twists a Demazure module by an anti-dominant character, then the resulting $B$-module can be filtered with successive quotients given by (isotypical) Demazure modules. This defines \emph{excellent filtrations} for $B$-modules which generalize \emph{good filtrations} for $G$-modules, in which successive quotients are Weyl modules. Tensor products of $G$-modules preserve good filtrations, but having excellent filtrations is not preserved by tensor products for $B$-modules \cite{vdK89}. In particular, $c_{\beta,\alpha}^{\gamma}$ are not, in general, nonnegative.

In this paper, we give explicit formulas for $c_{\beta,\alpha}^{\gamma}$ in the case when $\alpha = \nu$ is a partition of length $k$, in which case the corresponding Demazure character $\key_{\alpha}$ is the Schur polynomial $s_{\nu}(x_1,\ldots,x_k)$. Combinatorially, we describe $c_{\beta,\nu}^{\gamma}$ in
\begin{equation}
  \key_{\beta}(x_1,\ldots,x_n) s_{\nu}(x_1,\ldots,x_k) = \sum_{\gamma} c_{\beta,\nu}^{\gamma} \key_{\gamma}.
  \label{e:keyexps}
\end{equation}

Our new rule specializes to two previously known cases. Taking characters of Mathieu's excellent filtration \cite{Mat89} proves the product of a Demazure character with a Schur polynomial \emph{in at least as many variables}, that is, when $n\leq k$, expands nonnegatively into Demazure characters. Haglund, Luoto, Mason and van Willigenburg \cite{HLMvW11} give an explicit \emph{nonnegative} rule for the structure contants in this case. The author and Quijada \cite{AQ} give a \emph{signed} rule for \eqref{e:keyexps} in the Monk case, that is, when $\nu$ is degree $1$, \emph{in any number of variables}.

To realize positivity, we consider the geometric origins of Demazure modules. Recall $V(\lambda)$ may be constructed as dual to $H^0(G/B,\mathcal{L}_{\lambda})$, where $\mathcal{L}_{\lambda}$ is a certain line bundle over $G/B$ associated to $\lambda$. To each $w \in W$, we associate the \emph{Schubert cell} $C_w = B w B / B$, corresponding to the terms in the Bruhat decomposition, and the \emph{Schubert variety} $X_w$ which is its closure in the Zariski topology. Demazure \cite{Dem74a} considered the $B$-module $H^0(X_w,\mathcal{L}_{\lambda})$, which is dual to $V_w(\lambda)$.

Between the dual Weyl modules $H^0(G/B,\mathcal{L}_{\lambda})$ and the dual Demazure modules $H^0(X_w,\mathcal{L}_{\lambda})$ are the dual \emph{Schubert modules} $H^0(S, \mathcal{L}_{\lambda})$, where $S = \cup_{w\in \mathcal{I}} X_w$ is a union of Schubert varieties over a lower order ideal $\mathcal{I}\subset W$. Denote the Schubert modules by $V_\mathcal{I}(\lambda)$ and the \emph{Schubert characters} by $\key_{\Gamma} = \mathrm{ch}(V_{\mathcal{I}}(\lambda))$, where $\Gamma = \{ w\lambda \mid w \in \mathcal{I} \}$ is the lower order ideal in $\mathbb{Z}^n_{\ge 0}$. Polo \cite{Pol89} conjectured the following. 

\begin{conjecture*}[\cite{Pol89}]
  The tensor product of two dual Demazure modules has a \emph{Schubert filtration}, in which successive quotients are dual Schubert modules.
\end{conjecture*}

On the level of characters, this implies the product of Demazure characters is the sum of \emph{Demazure atoms} $\atom_{\gamma}$ whose indexing set forms a lower order ideal. Lascoux and Sch{\"u}tzenberger \cite{LS90} introduced Demazure atoms, under the name \emph{standard bases}, as the minimal non-intersecting pieces of Demazure characters. Thus Polo's conjecture implies the product of Demazure characters expands as a nonnegative sum of Demazure atoms, which was also directly conjectured by Pun \cite{Pun16}. We show our sum in \eqref{e:keyexps} over Demazure atoms can naturally be indexed by lower order ideals, proving both Pun's and Polo's conjectures for this case.

\begin{theorem*}
  We have explicit \emph{nonnegative} integers $a_{\beta,\nu}^{\gamma}$ and $a_{\beta,\nu}^{\Gamma}$ such that
\begin{equation}
  \key_{\beta}(x_1,\ldots,x_n) s_{\nu}(x_1,\ldots,x_k) = \sum_{\gamma} a_{\beta,\nu}^{\gamma} \atom_{\gamma}  = \sum_{\Gamma} a_{\beta,\nu}^{\Gamma} \key_{\Gamma},
  \label{e:atomexps}
\end{equation}  
\end{theorem*}

Our proofs are purely combinatorial, utilizing Kohnert's combinatorial model for Demazure characters \cite{Koh91}. Generalizing the Robinson--Schensted--Knuth insertion algorithm \cite{Rob38,Sch61,Knu70} on semistandard Young tableaux, we define an explicit insertion algorithm on Kohnert diagrams to give a bijection
\begin{equation}
  \KD(\beta) \times \mathrm{SSYT}(\nu) \xrightarrow{\sim}
  \bigsqcup_{\gamma} \left( \AKD(\gamma) \times \LAT(\gamma/\beta,\nu) \right) ,
  \label{e:k-rect}
\end{equation}
where $\KD(\beta)$ is the set of \emph{Kohnert diagrams} of shape $\beta$, $\AKD(\gamma)$ is the set of \emph{atomic Kohnert diagrams} of shape $\gamma$, and $\LAT(\gamma/\beta,\nu)$ is the set of \emph{lattice atomic tableaux} of \emph{skew shape} $\gamma/\beta$ and \emph{weight} $\nu$. The coefficient $a_{\beta,\nu}^{\gamma}$ in \eqref{e:atomexps} is the cardinality of $\LAT(\gamma/\beta,\nu)$, and $a_{\beta,\nu}^{\Gamma}$ is the cardinality of $\LAT(\gamma/\beta,\nu)$ for any generator $\gamma \in \Gamma$. We derive our signed formula for $c_{\beta,\nu}^{\gamma}$ in \eqref{e:keyexps} from this via inclusion--exclusion.

The combinatorics developed in this paper has been extended to an analogous conjecture for Schubert times Schur structure constants \cite{A-pnas}, with a proof announced by the author and Bergeron generalizing this insertion algorithm \cite{AB}.

%
\section{Kohnert diagrams}
%
\label{sec:kohnert}

In this section, we give combinatorial models for polynomial bases using \newword{diagrams}, finite collections of \emph{cells} in the first quadrant of $\mathbb{Z}\times\mathbb{Z}$, and we introduce a new combinatorially defined basis of \emph{pinned polynomials}.

\subsection{Kohnert's rule}
\label{sec:diagrams}

Demazure \cite{Dem74} generalized the Weyl character formula to Demazure modules, with proof gaps noted by Kac and completed by Joseph \cite{Jos85} using excellent filtrations. Kohnert's elegant combinatorial model for Demazure characters \cite{Koh91} uses the following operation on diagrams.

\begin{definition}[\cite{Koh91}]
  A \newword{Kohnert move} on a diagram selects the rightmost cell of a given row and moves the cell down within its column to the first available position below, if it exists, jumping over other cells in its way as needed. 
  \label{def:kohnertmove}
\end{definition}

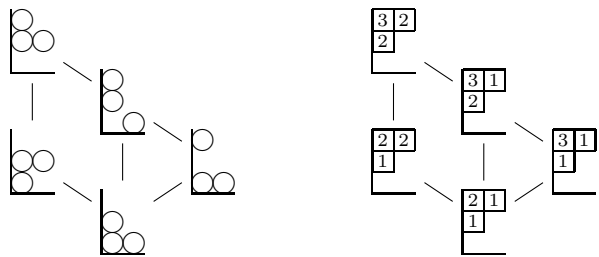
\begin{figure}[ht]
  \begin{center}
    \begin{tikzpicture}[xscale=1.2,yscale=0.8]
      \node at (0,5) (A) {$\cirtab{ ~ \\ ~ & ~ \\ \\\hline}$};
      \node at (1,4) (B) {$\cirtab{ ~ \\ ~ \\ & ~ \\\hline}$};
      \node at (2,3) (C) {$\cirtab{ ~ \\ \\ ~ & ~ \\\hline}$};
      \node at (1,2) (E) {$\cirtab{ \\ ~ \\ ~ & ~ \\\hline}$};
      \node at (0,3) (H) {$\cirtab{ \\ ~ & ~ \\ ~ \\\hline}$};
      \draw[thin] (A) -- (H) ;
      \draw[thin] (A) -- (B) ;
      \draw[thin] (H) -- (E) ;
      \draw[thin] (B) -- (E) ;
      \draw[thin] (B) -- (C) ;
      \draw[thin] (C) -- (E) ;
      \node at (4,5) (A2) {$\vline\tableau{3 & 2 \\ 2 \\ & \\\hline}$};
      \node at (5,4) (B2) {$\vline\tableau{3 & 1 \\ 2 \\ & \\\hline}$};
      \node at (6,3) (C2) {$\vline\tableau{3 & 1 \\ 1 \\ & \\\hline}$};
      \node at (5,2) (E2) {$\vline\tableau{2 & 1 \\ 1 \\ & \\\hline}$};
      \node at (4,3) (H2) {$\vline\tableau{2 & 2 \\ 1 \\ & \\\hline}$};
      \draw[thin] (A2) -- (H2) ;
      \draw[thin] (A2) -- (B2) ;
      \draw[thin] (H2) -- (E2) ;
      \draw[thin] (B2) -- (E2) ;
      \draw[thin] (B2) -- (C2) ;
      \draw[thin] (C2) -- (E2) ;
    \end{tikzpicture}
    \caption{\label{fig:kohnert}The set $\KD(0,2,1)$ of Kohnert diagrams for $(0,2,1)$ and their images under the injection with the set $\SSRT(0,1,2)$ of semistandard reverse tableaux of partition shape $(0,1,2)$.}
  \end{center}
\end{figure}

The set of \newword{Kohnert diagrams} for $D$, denoted by $\KD(D)$, consists of all diagrams obtainable from $D$ by some (possibly empty) sequence of Kohnert moves; see Fig.~\ref{fig:kohnert}. The \newword{composition diagram} $\D(\alpha)$ for a weak composition $\alpha$ has $\alpha_i$ cells, left-justified in row $i$. We compress notation by $\KD(\alpha) = \KD(\D(\alpha))$.

\begin{definition}[\cite{Koh91}]
  The \newword{Demazure character} $\key_{\alpha}$ for $\alpha \in \mathbb{N}^m$ is 
  \begin{equation}
    \key_{\alpha}(x_1,\ldots,x_m) = \sum_{T \in \KD(\alpha)} x_1^{\wt(T)_1} \cdots x_m^{\wt(T)_m},
    \label{e:key-d}
  \end{equation}
  where $\wt(T)_r$ is the number of cells in row $r$ of $T$.
  \label{def:kohnert}
\end{definition}

Not all diagrams can arise as Kohnert diagrams for some weak composition.

\begin{definition}
  A diagram $T$ is \newword{rectified} if $T \in \KD(\alpha)$ for some $\alpha$.
  \label{def:generic}
\end{definition}

A \newword{partition} is a weakly \emph{increasing} weak composition. A \newword{semistandard reverse tableau of partition shape $\lambda$}, denoted by $T\in\SSRT(\lambda)$, is a filling of the cells of $\D(\lambda)$ with integers $1 \le i \le \ell(\lambda)$ such that entries weakly decrease left to right within rows and strictly decrease top to bottom within columns; see Fig.~\ref{fig:kohnert}.

\begin{definition}
  The \newword{Schur polynomial} indexed by the partition $\lambda\in\mathbb{N}^k$ is
  \begin{equation}
    s_{\lambda}(x_1,\ldots,x_k) = \sum_{T \in \SSRT(\lambda)} x_1^{\wt(T)_1} \cdots x_k^{\wt(T)_k},
    \label{e:schur}
  \end{equation}
  where $k = \ell(\lambda)$ and $\wt(T)_i$ is the number of entries of $T$ equal to $i$.
  \label{def:schur}
\end{definition}

To relate the diagram model with the familiar paradigm of tableaux, Assaf and Searles \cite[Thm~4.6]{AS18} prove the following, rephrased here for reverse tableaux. 

\begin{proposition}[\cite{AS18}]
  For $\alpha\in\mathbb{N}^m$ and $\sort(\alpha)$ the unique partition in the $\mathcal{S}_{m}$ orbit of $\alpha$, the map $\varphi:\KD(\alpha) \rightarrow \SSRT(\sort(\alpha))$ that assigns entry $i$ to each cell in row $i$ and lifts the cells within their columns is a weight-preserving, injective map that is surjective if and only if $\alpha = \sort(\alpha)$.
  \label{prop:raise}
\end{proposition}

In particular, every Schur polynomial is a Demazure character. Moreover, under this injection, our new constructions align with classical operations on tableaux.

\subsection{Threading algorithm}
\label{sec:threading}

\begin{figure}[ht]
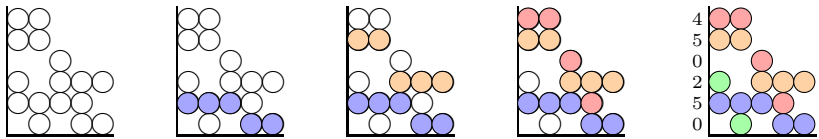

  \begin{displaymath}
    \arraycolsep=1.5\cellsize
    \begin{array}{ccccc}
      \cirtab{
        {~} & {~} \\
        {~} & {~} \\
        & & {~} \\
        {~} & & {~} & {~} & {~} \\
        {~} & {~} & {~} & {~} \\
        & {~} & & {~} & {~} \\\hline} &
      \cirtab{
        {~} & {~} \\
        {~} & {~} \\
        & & {~} \\
        {~} & & {~} & {~} & {~} \\
        \cball{\ctwo}{~} & \cball{\ctwo}{~} & \cball{\ctwo}{~} & {~} \\
        & {~} & & \cball{\ctwo}{~} & \cball{\ctwo}{~} \\\hline} &
      \cirtab{
        {~} & {~} \\
        \cball{\cfiv}{~} & \cball{\cfiv}{~} \\
        & & {~} \\
        {~} & & \cball{\cfiv}{~} & \cball{\cfiv}{~} & \cball{\cfiv}{~} \\
        \cball{\ctwo}{~} & \cball{\ctwo}{~} & \cball{\ctwo}{~} & {~} \\
        & {~} & & \cball{\ctwo}{~} & \cball{\ctwo}{~} \\\hline} &
      \cirtab{
        \cball{\csix}{~} & \cball{\csix}{~} \\
        \cball{\cfiv}{~} & \cball{\cfiv}{~} \\
        & & \cball{\csix}{~} \\
        {~} & & \cball{\cfiv}{~} & \cball{\cfiv}{~} & \cball{\cfiv}{~} \\
        \cball{\ctwo}{~} & \cball{\ctwo}{~} & \cball{\ctwo}{~} & \cball{\csix}{~} \\
        & {~} & & \cball{\ctwo}{~} & \cball{\ctwo}{~} \\\hline} &
      \nulltab{4 \\ 5 \\ 0 \\ 2 \\ 5 \\ 0}\hss\vline\nulltab{
        \cball{\csix}{~} & \cball{\csix}{~} \\
        \cball{\cfiv}{~} & \cball{\cfiv}{~} \\
        & & \cball{\csix}{~} \\
        \cball{\cthr}{~} & & \cball{\cfiv}{~} & \cball{\cfiv}{~} & \cball{\cfiv}{~} \\
        \cball{\ctwo}{~} & \cball{\ctwo}{~} & \cball{\ctwo}{~} & \cball{\csix}{~} \\
        & \cball{\cthr}{~} & & \cball{\ctwo}{~} & \cball{\ctwo}{~} \\\hline} 
    \end{array}
  \end{displaymath}
    \caption{\label{fig:threading}Constructing the thread decomposition of a rectified diagram.}
\end{figure}

Demazure atoms were introduced by Lascoux and Sch{\"u}tzenberger \cite{LS90} as a refinement of Demazure characters. Mason \cite{Mas09} gave a tableau model for Demazure atoms, and Searles \cite{Sea20} gave a diagram model based on the \emph{threading algorithm} of Assaf and Searles \cite[Def~3.5]{AS18}; see Fig.~\ref{fig:threading}.

\begin{definition}[\cite{AS18}]
  The \newword{thread decomposition} of a diagram partitions cells into \newword{threads} as follows. Begin a new thread with the rightmost then lowest unthreaded cell. After threading a cell in column $c+1$, thread the lowest unthreaded cell in column $c$ weakly above it, and continue left until column $1$ or until no cell can be threaded. Continue threading until all cells are threaded.
  \label{def:threading}
\end{definition}

Assaf and Searles prove every thread of a diagram $T$ ends in the first column if and only if $T$ is rectified \cite[Lemma~2.2]{AS18}, showing the following is well-defined.

\begin{definition}[\cite{AS18}]
  The \newword{thread weight} of a rectified diagram $T$, denoted by $\thread(T)$, has $\thread(T)_r$ equal to the number of cells in the thread in column $1$, row $r$.
\end{definition}

Assaf and Searles \cite{AS18} use threads to partition the set $\KD(\alpha)$ corresponding to the decomposition of a Demazure character into quasi-key polynomials. Searles \cite{Sea20} refines this to reflect the decomposition of a Demazure character into atoms.

\begin{definition}
  A rectified diagram $T$ is an \newword{atomic Kohnert diagram for $\alpha$} if $\thread(T) = \alpha$. Denote the set of atomic Kohnert diagrams for $\alpha$ by $\AKD(\alpha)$.
  \label{def:AKD}
\end{definition}

As this will be important for deriving formulas from the combinatorics to follow, notice for $\alpha \neq \beta$, the sets $\AKD(\alpha)$ and $\AKD(\beta)$ are \emph{disjoint}.

Combining \cite[Thm~3.7]{AS18} and \cite[Thm~3.6]{Sea20} gives the following.

\begin{definition}[\cite{AS18,Sea20}]
  The \newword{Demazure atom} $\atom_{\alpha}$ for $\alpha \in \mathbb{N}^m$ is 
  \begin{equation}
    \atom_{\alpha}(x_1,\ldots,x_m) = \sum_{T \in \AKD(\alpha)} x_1^{\wt(T)_1} \cdots x_m^{\wt(T)_m} 
    \label{e:atom-d}
  \end{equation}
  where $\wt(T)_r$ is the number of cells in row $r$ of $T$.
  \label{def:atom-d}
\end{definition}

Following \cite{AS18}, a \newword{left swap} on a weak composition $\alpha$ exchanges two parts $\alpha_r < \alpha_s$ with $r<s$. Write $\alpha \lswap \beta$ whenever $\alpha$ is obtainable via some (possibly empty) sequence of left swaps on $\beta$. The following is implicit in \cite[Thm~3.7]{AS18}.

\begin{lemma}[\cite{AS18}]
  For $T$ rectified, $T\in\KD(\alpha)$ if and only if $\thread(T) \lswap \alpha$.
  \label{lem:lswap}
\end{lemma}

Thus we can refine Kohnert diagrams into their atomic subsets by
\begin{equation}
  \KD(\beta) = \bigsqcup_{\alpha \lswap \beta} \AKD(\alpha) .
  \label{e:KD2AKD}
\end{equation}
Taking generating polynomials gives the familiar decomposition
\begin{equation}
  \key_{\beta} = \sum_{\alpha \lswap \beta} \atom_{\alpha}.
  \label{e:key2atom}
\end{equation}
That is, Demazure characters are sums over Bruhat intervals of Demazure atoms.

\begin{definition}
  The \newword{Schubert character} for $\Gamma\subset\mathbb{N}^m$ a lower order ideal  is
  \begin{equation}
    \key_{\Gamma} = \sum_{\alpha\in \Gamma} \atom_{\alpha}.
    \label{e:schubchar}
  \end{equation}
\end{definition}

Then Polo's conjecture \cite{Pol89} for characters can be stated as follows.

\begin{conjecture}[\cite{Pol89}]
  There exist nonnegative integers $a_{\beta,\alpha}^{\Gamma}$ such that 
  \begin{equation}
    \key_{\beta} \key_{\alpha} = \sum_{\Gamma} a_{\beta,\alpha}^{\Gamma} \key_{\Gamma}.
    \label{e:schubchar-general}
  \end{equation}
  \label{conj:schubchar}
\end{conjecture}

Since Schubert characters over determine a basis, there is no computational test for Conjecture~\ref{conj:schubchar}. However, by \eqref{e:schubchar}, the right hand side of \eqref{e:schubchar-general} expands nonnegatively into Demazure atoms, which are a basis. This motivates the following weaker conjecture \cite[Conj.~1]{Pun16} for which Pun gives a proof for compositions $\beta,\alpha$ of length at most $3$ with at most $2$ nonzero parts.

\begin{conjecture}[\cite{Pun16}]
  There exist nonnegative integers $a_{\beta,\alpha}^{\gamma}$ such that 
  \begin{equation}
    \key_{\beta} \key_{\alpha} = \sum_{\gamma} a_{\beta,\alpha}^{\gamma} \atom_{\gamma}.
    \label{e:atom-general}
  \end{equation}
  \label{conj:atom}
\end{conjecture}

In Theorems~\ref{thm:atom} and \ref{thm:schubchar}, we give nonnegative formulas for $a_{\beta,\alpha}^{\gamma}$ and $a_{\beta,\alpha}^{\Gamma}$, respectively, when $\alpha$ is a partition, thereby resolving both conjectures in this case.

\subsection{Labeled diagrams}
\label{sec:labels}

Assaf and Searles give a canonical labeling \cite[Def~2.3]{AS18} of cells of a rectified diagram by a composition $\alpha$. 

\begin{definition}[\cite{AS18}]
  Given a diagram, a \newword{semi-proper labeling of shape $\alpha$} assigns positive integers to the cells of the diagram such that 
  \begin{enumerate}    
  \item \label{i:shape} [\newword{strict}] column $c$ consists of distinct entries $\{ r \mid \alpha_r \ge c\}$;
  \item \label{i:flag} [\newword{flagged}] each entry in row $r$ is at least $r$;
  \item \label{i:descend} [\newword{descending}] cells with entry $r$ weakly descend from left to right.
  \end{enumerate}
  A labeling is \newword{proper} if in addition
  \begin{enumerate}    
    \setcounter{enumi}{3}
  \item \label{i:invert} [\newword{minimal}] if $r<s$ appear in a column with $r$ above $s$, then there is an $r$ in the column immediately to the right of and strictly above $s$.
  \end{enumerate}
  \label{def:KT}
\end{definition}

\begin{figure}[ht]
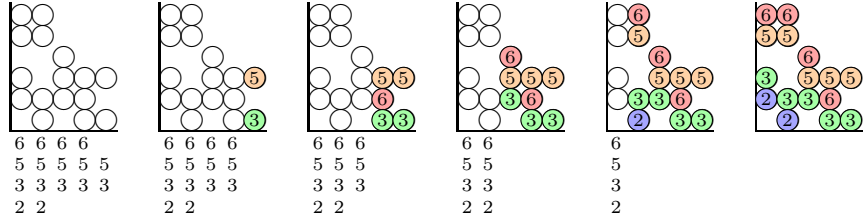

  \begin{displaymath}
    \arraycolsep=\cellsize
    \begin{array}{llllll}
      \cirtab{
        {~} & {~} \\
        {~} & {~} \\
        & & {~} \\
        {~} & & {~} & {~} & {~} \\
        {~} & {~} & {~} & {~} \\
        & {~} & & {~} & {~} \\\hline} &
      \cirtab{
        {~} & {~} \\
        {~} & {~} \\
        & & {~} \\
        {~} & & {~} & {~} & \cball{\cfiv}{5} \\
        {~} & {~} & {~} & {~} \\
        & {~} & & {~} & \cball{\cthr}{3} \\\hline} &
      \cirtab{
        {~} & {~} \\
        {~} & {~} \\
        & & {~} \\
        {~} & & {~} & \cball{\cfiv}{5} & \cball{\cfiv}{5} \\
        {~} & {~} & {~} & \cball{\csix}{6} \\
        & {~} & & \cball{\cthr}{3} & \cball{\cthr}{3} \\\hline} &
      \cirtab{
        {~} & {~} \\
        {~} & {~} \\
        & & \cball{\csix}{6} \\
        {~} & & \cball{\cfiv}{5} & \cball{\cfiv}{5} & \cball{\cfiv}{5} \\
        {~} & {~} & \cball{\cthr}{3} & \cball{\csix}{6} \\
        & {~} & & \cball{\cthr}{3} & \cball{\cthr}{3} \\\hline} &
      \cirtab{
        {~} & \cball{\csix}{6} \\
        {~} & \cball{\cfiv}{5} \\
        & & \cball{\csix}{6} \\
        {~} & & \cball{\cfiv}{5} & \cball{\cfiv}{5} & \cball{\cfiv}{5} \\
        {~} & \cball{\cthr}{3} & \cball{\cthr}{3} & \cball{\csix}{6} \\
        & \cball{\ctwo}{2} & & \cball{\cthr}{3} & \cball{\cthr}{3} \\\hline} &
      \cirtab{
        \cball{\csix}{6} & \cball{\csix}{6} \\
        \cball{\cfiv}{5} & \cball{\cfiv}{5} \\
        & & \cball{\csix}{6} \\
        \cball{\cthr}{3} & & \cball{\cfiv}{5} & \cball{\cfiv}{5} & \cball{\cfiv}{5} \\
        \cball{\ctwo}{2} & \cball{\cthr}{3} & \cball{\cthr}{3} & \cball{\csix}{6} \\
        & \cball{\ctwo}{2} & & \cball{\cthr}{3} & \cball{\cthr}{3} \\\hline}  \\
      \nulltab{6 & 6 & 6 & 6 \\ 5 & 5 & 5 & 5 & 5 \\ 3 & 3 & 3 & 3 & 3 \\ 2 & 2} &      
      \nulltab{6 & 6 & 6 & 6 \\ 5 & 5 & 5 & 5 \\ 3 & 3 & 3 & 3 \\ 2 & 2} &      
      \nulltab{6 & 6 & 6  \\ 5 & 5 & 5  \\ 3 & 3 & 3  \\ 2 & 2} &      
      \nulltab{6 & 6   \\ 5 & 5   \\ 3 & 3  \\ 2 & 2} &      
      \nulltab{6 \\ 5 \\ 3 \\ 2 } &      
    \end{array}
  \end{displaymath}
    \caption{\label{fig:labeling}Constructing the proper labeling of a rectified diagram for the composition $\alpha = (0,2,5,0,5,4)$.}
\end{figure}

Assaf and Searles \cite[Def~2.5]{AS18} construct the proper labeling as follows; see Fig.~\ref{fig:labeling}.

\begin{definition}[\cite{AS18}]
  The \newword{proper labeling of $T$ with respect to $\alpha$}, denoted by $\Label_{\alpha}(T)$, is defined by: once all columns right of $c$ have been labeled, bijectively assign labels $\{i \mid a_i \geq c\}$ to cells of column $c$ from bottom to top by choosing the smallest unused label $i$ such that the $i$ in column $c+1$, if it exists, is weakly lower.
  \label{def:kohnert-label}
\end{definition}

Combining \cite[Thm~2.8]{AS18} and \cite[Lemma~5.2.1]{Ass22-KC} gives the following.

\begin{lemma}[\cite{AS18,Ass22-KC}]
  For a diagram $T$, the following are equivalent
  \begin{enumerate}[label=(\roman*)]
  \item $T\in\KD(\alpha)$;
  \item the proper labeling $\Label_{\alpha}(T)$ is well-defined and flagged;
  \item $T$ has a semi-proper labeling of shape $\alpha$.
  \end{enumerate}
  \label{lem:matching}
\end{lemma}

\begin{figure}[ht]
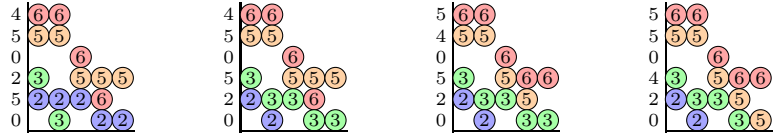

  \begin{displaymath}
    \arraycolsep=2\cellsize
    \begin{array}{cccc}
      \nulltab{4 \\ 5 \\ 0 \\ 2 \\ 5 \\ 0}\hss\vline\nulltab{
        \cball{\csix}{6} & \cball{\csix}{6} \\
        \cball{\cfiv}{5} & \cball{\cfiv}{5} \\ & & \cball{\csix}{6} \\
        \cball{\cthr}{3} & & \cball{\cfiv}{5} & \cball{\cfiv}{5} & \cball{\cfiv}{5} \\
        \cball{\ctwo}{2} & \cball{\ctwo}{2} & \cball{\ctwo}{2} & \cball{\csix}{6} \\
        & \cball{\cthr}{3} & & \cball{\ctwo}{2} & \cball{\ctwo}{2} \\\hline} &
      \nulltab{4 \\ 5 \\ 0 \\ 5 \\ 2 \\ 0}\hss\vline\nulltab{
        \cball{\csix}{6} & \cball{\csix}{6} \\
        \cball{\cfiv}{5} & \cball{\cfiv}{5} \\ & & \cball{\csix}{6} \\
        \cball{\cthr}{3} & & \cball{\cfiv}{5} & \cball{\cfiv}{5} & \cball{\cfiv}{5} \\
        \cball{\ctwo}{2} & \cball{\cthr}{3} & \cball{\cthr}{3} & \cball{\csix}{6} \\
        & \cball{\ctwo}{2} & & \cball{\cthr}{3} & \cball{\cthr}{3} \\\hline} &
      \nulltab{5 \\ 4 \\ 0 \\ 5 \\ 2 \\ 0}\hss\vline\nulltab{
        \cball{\csix}{6} & \cball{\csix}{6} \\
        \cball{\cfiv}{5} & \cball{\cfiv}{5} \\ & & \cball{\csix}{6} \\
        \cball{\cthr}{3} & & \cball{\cfiv}{5} & \cball{\csix}{6} & \cball{\csix}{6} \\
        \cball{\ctwo}{2} & \cball{\cthr}{3} & \cball{\cthr}{3} & \cball{\cfiv}{5} \\
        & \cball{\ctwo}{2} & & \cball{\cthr}{3} & \cball{\cthr}{3} \\\hline} &
      \nulltab{5 \\ 5 \\ 0 \\ 4 \\ 2 \\ 0}\hss\vline\nulltab{
        \cball{\csix}{6} & \cball{\csix}{6} \\
        \cball{\cfiv}{5} & \cball{\cfiv}{5} \\ & & \cball{\csix}{6} \\
        \cball{\cthr}{3} & & \cball{\cfiv}{5} & \cball{\csix}{6} & \cball{\csix}{6} \\
        \cball{\ctwo}{2} & \cball{\cthr}{3} & \cball{\cthr}{3} & \cball{\cfiv}{5} \\
        & \cball{\ctwo}{2} & & \cball{\cthr}{3} & \cball{\cfiv}{5} \\\hline} 
    \end{array}
  \end{displaymath}
    \caption{\label{fig:pinned-labeling}Pinned, proper labelings of a rectified diagram $T$; the leftmost is the atomic labeling showing the thread decomposition.}
\end{figure}

The \newword{atomic labeling of $T$}, denoted by $\Label_{\thread}(T)$, in which each cell along a thread ending in row $r$ is labeled $r$, is proper and \emph{pinned}; see Fig.~\ref{fig:pinned-labeling}.

\begin{definition}
  A labeling is \newword{pinned} if each entry in column $1$ equals its row.  Given $\alpha$, denote the set of \newword{pinned Kohnert diagrams of shape $\alpha$} with pinned, proper labels by $\PKD(\alpha) = \{T\in\KD(\alpha)\mid \Label_{\alpha}(T) \text{ is pinned}\}$.
  \label{def:PKT}
\end{definition}

A \newword{pinned left swap} on $\beta$ exchanges two \emph{nonzero} parts $0 < \beta_i < \beta_j$ for $i<j$. Write $\alpha \lswappin \beta$ whenever $\alpha$ results from some sequence of pinned left swaps on $\beta$.

\begin{lemma}
  We have $\AKD(\alpha) \subseteq \PKD(\beta)$ if and only if $\alpha \lswappin \beta$.
  \label{lem:lswappin}
\end{lemma}

\begin{proof}
  If $\alpha \lswappin \beta$, then $\alpha\lswap\beta$ so, by Lemma~\ref{lem:lswap}, $\KD(\alpha)\subseteq\KD(\beta)$. Thus $T\in\AKD(\alpha)\subseteq\KD(\beta)$ and since $\alpha \lswappin \beta$, $T$ is also pinned for $\beta$.

  Conversely, suppose $\AKD(\alpha) \subseteq \PKD(\beta)$. Then $\alpha_i = 0$ if and only if $\beta_i=0$. Since $\D(\alpha)\in\PKD(\beta)\subseteq\KD(\beta)$, by Lemma~\ref{lem:lswap}, $\alpha=\thread(\D(\alpha)) \lswap\beta$ as well.
\end{proof}

The generating polynomial for pinned Kohnert diagrams provides a useful intermediate basis between Demazure characters and Demazure atoms. 

\begin{definition}
  The \newword{pinned polynomial} $\pinp_{\alpha}$ for $\alpha \in \mathbb{N}^m$ is 
  \begin{equation}
    \pinp_{\alpha}(x_1,\ldots,x_m) = \sum_{T \in \PKD(\alpha)} x_1^{\wt(T)_1} \cdots x_m^{\wt(T)_m}.
    \label{e:pinp}
  \end{equation}
  \label{def:pinp}
\end{definition}


\begin{proposition}
  The set $\{\pinp_{\beta}\}_{\beta\in\mathbb{N}^m}$ is a $\mathbb{Z}$-basis for $\mathbb{Z}[x_1,\ldots,x_m]$. Moreover, for $\alpha,\beta,\gamma$ compositions of length $m$, we have nonnegative expansions
  \begin{equation}
    \key_{\gamma} = \sum_{\substack{\beta\lswap\gamma \\ \beta\lswappin\alpha\lswap\gamma \Rightarrow \alpha=\beta}}  \pinp_{\beta}
    \hspace{1em} \text{and} \hspace{1em}
    \pinp_{\beta} = \sum_{\alpha \lswappin \beta} \atom_{\alpha}.
    \label{e:2pin}
  \end{equation}
  \label{prop:basis}
\end{proposition}

\begin{proof}
  The right side of \eqref{e:2pin} follows from Lemma~\ref{lem:lswappin}. The left expansion follows from this and \eqref{e:KD2AKD} together with the observation that for any $\alpha$, there exists a unique maximal $\beta$ for which $\alpha\lswappin\beta\lswap\gamma$. In particular, the pinned polynomials are upper uni-triangular with respect to Demazure atoms, which are a $\mathbb{Z}$-basis. 
\end{proof}

The product of a Demazure atom with a Schur polynomial does not always expand nonnegatively into Demazure atoms. The motivation for introducing this new basis is that the product of a pinned polynomial and a Schur polynomial does. 

%
\section{Insertion algorithm}
%
\label{sec:insertion}

In this section, we define our algorithm for inserting a cell into labeled Kohnert diagrams where we bound the insertion by a parameter $k$. We show this is equivalent to RSK via the injection $\varphi$ when $k$ is sufficiently large.

\subsection{Label exchanges}
\label{sec:exchange}

We restrict our inserted cell by a parameter $k$ and take into account the pinned labeling of the diagram into which it is inserted. 

\begin{definition}
  For $r<s$ in a semi-proper labeling, $s$ \newword{touches} $r$ at column $c$ if there is an $s$ in column $c+1$ and it lies weakly below $r$ in column $c$. We say $s$ \newword{crosses} $r$ if there is an $r$ above $s$ in column $c+1$; otherwise $s$ \newword{abuts} $r$.
  \label{def:cross}
\end{definition}

\begin{figure}[ht]
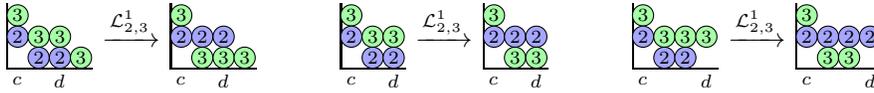

  \begin{displaymath}
    \arraycolsep=2pt
    \begin{array}{ccccccccccc}
      \vline\nulltab{
        \cball{\cthr}{3} & \\
        \cball{\ctwo}{2} & \cball{\cthr}{3} & \cball{\cthr}{3} \\
        & \cball{\ctwo}{2} & \cball{\ctwo}{2} & \cball{\cthr}{3} \\\hline}
      & \raisebox{-\cellsize}{$\xrightarrow{\Label_{2,3}^{1}}$} &
      \vline\nulltab{
        \cball{\cthr}{3} & \\
        \cball{\ctwo}{2} & \cball{\ctwo}{2} & \cball{\ctwo}{2} \\
        & \cball{\cthr}{3} & \cball{\cthr}{3} & \cball{\cthr}{3} \\\hline} 
      & \hspace{3\cellsize} &
      \vline\nulltab{
        \cball{\cthr}{3} & \\
        \cball{\ctwo}{2} & \cball{\cthr}{3} & \cball{\cthr}{3}  \\
        & \cball{\ctwo}{2} & \cball{\ctwo}{2} \\\hline}
      & \raisebox{-\cellsize}{$\xrightarrow{\Label_{2,3}^{1}}$} &
      \vline\nulltab{
        \cball{\cthr}{3} & \\
        \cball{\ctwo}{2} & \cball{\ctwo}{2} & \cball{\ctwo}{2} \\
        & \cball{\cthr}{3} & \cball{\cthr}{3} \\\hline}
      & \hspace{3\cellsize} &
      \vline\nulltab{
        \cball{\cthr}{3} & \\
        \cball{\ctwo}{2} & \cball{\cthr}{3} & \cball{\cthr}{3} & \cball{\cthr}{3} \\
        & \cball{\ctwo}{2} & \cball{\ctwo}{2} \\\hline}
      & \raisebox{-\cellsize}{$\xrightarrow{\Label_{2,3}^{1}}$} &
      \vline\nulltab{
        \cball{\cthr}{3} & \\
        \cball{\ctwo}{2} & \cball{\ctwo}{2} & \cball{\ctwo}{2} & \cball{\ctwo}{2} \\
        & \cball{\cthr}{3} & \cball{\cthr}{3} \\\hline} \\
      \nulltab{ c & & d & } & &
      \nulltab{ c & & d & } & &
      \nulltab{ c & & d } & &
      \nulltab{ c & & d } & &
      \nulltab{ c & & & d} & &
      \nulltab{ c & & & d} 
    \end{array}
  \end{displaymath}
    \caption{\label{fig:touch}Labels that touch at column $1$ and their exchanges.}
\end{figure}

The following procedure turns a proper labeling into a semi-proper labeling with smaller labels raised higher and the shape related by a pinned left swap.

\begin{definition}
  If $s$ abuts $r$ in column $c$ of a semi-properly labeled diagram, the \newword{exchange labeling} for $r,s$ at column $c$, denoted by $\Label_{r,s}^{c}$, is obtained by simultaneously changing $r$ to $s$ and $s$ to $r$ in columns $c+1,\ldots,d$, where $d$ is the leftmost column right of $c$ in which $s$ touches $r$, if it exists, or else $d = \alpha_s$.
\label{def:exchange-lab}
\end{definition}

\begin{lemma}
  Exchange labelings of semi-proper labelings are semi-proper.
  \label{lem:exchange}
\end{lemma}

\begin{proof}
  For condition \eqref{i:shape}, if $s$ touches $r$ after column $c$ or if $\alpha_s \le \alpha_r$, then the set of labels within each column is unchanged, and so the weight of the labeling is $\alpha$. Otherwise, labels $s$ terminate at $\min(\alpha_r,\alpha_s)$ and labels $r$ extend to $\max(\alpha_r,\alpha_s)$, changing the weight to $t_{r,s} \alpha$ if $\alpha_r<\alpha_s$, but still maintaining \eqref{i:shape}. 

  For condition \eqref{i:flag}, the leftmost affected labels occur weakly below a row $q$ occupied by a label $r$, and so $q<r$. By \eqref{i:descend} for the original labeling, every affected label occurs weakly below row $q$, and so \eqref{i:flag} holds for the exchange labeling.

  Finally, condition \eqref{i:descend} holds in columns $b \le c$ and $b > d$ since those columns are unaffected, and in columns $c+1 < b \le d$ since the labels are swapped in those columns. Thus we need only check at columns $c,c+1$ and $d,d+1$. From the definition of touching, both $r,s$ in column $c+1$ (if they exist) lie weakly below both $r,s$ in column $c$ as required. Similarly, if $s$ touches $r$ at column $d$, then both $r,s$ in column $d+1$ (if they exist) lie weakly below both $r,s$ in column $d$ as required. If $s$ never again touches $r$, then $s$ terminates at $\min(\alpha_r,\alpha_s)\le d$, leaving nothing to show for $s$. In this case, if $\alpha_r \le \alpha_s=d$, there is also nothing to show for $r$, and if $\alpha_r > \alpha_s$, then since $s,r$ never again touched, $s$ lies above $r$ in column $d$, and so above $r$ in column $d+1$ by condition \eqref{i:descend} for the original labeling. Thus \eqref{i:descend} holds for the exchange labeling as well, and so it is semi-proper.
\end{proof}

We now prove that semi-proper labelings are labelings which can be obtained from the proper labeling by exchange labelings .

\begin{theorem}
  For $T\in\PKD(\alpha)$, every semi-proper, pinned labeling of $T$ of shape $\alpha'\lswappin\alpha$ can be obtained from $\Label_{\alpha}(T)$ via exchange labelings.
  \label{thm:semiproper}
\end{theorem}

\begin{proof}
  Consider the cells $x_1,\ldots,x_N$ of $T$ taken column by column from top to bottom, left to right. Given two distinct labelings, set $\delta(\Label, \Label') = d$ if $\Label(x_i) = \Label'(x_i)$ for all $i<d$ and $\Label(x_d) \neq \Label'(x_d)$; that is, $x_d$ is the first cell at which the labelings differ. Extend this to a set $L$ of labelings by
  \[ \delta(\Label,L) = \max \{ \delta(\Label, \Label') \mid \Label'\in L \}. \]

  Let $L$ denote the subset of semi-proper, pinned labelings which can be obtained from the proper labeling via exchange labelings. Suppose, for contradiction, there exist some semi-proper, pinned labeling not in $L$. Take $\Label \not\in L$, and choose $\Label'\in L$ such that $\delta = \delta(\Label,\Label') = \delta(\Label,L)$. Let $c+1$ be the column of $x_{\delta}$. Since all pinned labelings must agree in column $1$, we must have $c\ge 1$. 

  Suppose $r = \Label(x_{\delta}) < \Label'(x_{\delta}) = s$. Since $\Label$ is descending, the cell with label $r$ in column $c$ lies weakly above $x_{\delta}$ in both labelings (which agree prior to $x_{\delta}$). Since $\Label$ is strict, there there must be a cell $y$ in column $c+1$ below $x_{\delta}$ with $\Label'(y)=r$. In particular, $s$ abuts $r$ in column $c$. Therefore, the exchange $\Label_{r,s}^{c}$ on $\Label'$ is well-defined, and so lies in $L$. Moreover, $\delta(\Label,\Label_{r,s}^{c}) > \delta$, a contradiction. The case when $r>s$ is symmetric, since both labelings are semi-proper.
\end{proof}

\begin{figure}[ht]
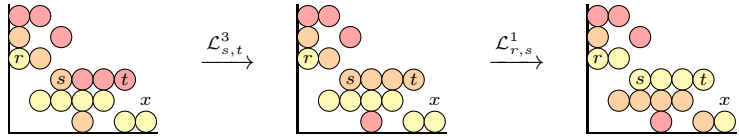

  \begin{displaymath}
    \arraycolsep=\cellsize
    \begin{array}{ccccc}
      \vline\nulltab{
        \cball{\csix}{} & \cball{\csix}{} \\
        \cball{\cfiv}{} & & \cball{\csix}{} \\
        \cball{\cfou}{r} & \cball{\cfiv}{} \\
        & & \cball{\cfiv}{s} & \cball{\csix}{} & \cball{\csix}{} & \cball{\csix}{t} \\
        & \cball{\cfou}{} & \cball{\cfou}{} & \cball{\cfou}{} & \cball{\cfou}{} & & x \\
        & & & \cball{\cfiv}{} & & \cball{\cfou}{} & \cball{\cfou}{} \\\hline }
      & \raisebox{-2\cellsize}{$\xrightarrow{\Label_{s,t}^{3}}$} &
      \vline\nulltab{
        \cball{\csix}{} & \cball{\csix}{} \\
        \cball{\cfiv}{} & & \cball{\csix}{} \\
        \cball{\cfou}{r} & \cball{\cfiv}{} \\
        & & \cball{\cfiv}{s} & \cball{\cfiv}{} & \cball{\cfiv}{} & \cball{\cfiv}{t} \\
        & \cball{\cfou}{} & \cball{\cfou}{} & \cball{\cfou}{} & \cball{\cfou}{} & & x \\
        & & & \cball{\csix}{} & & \cball{\cfou}{} & \cball{\cfou}{} \\\hline }
      & \raisebox{-2\cellsize}{$\xrightarrow{\Label_{r,s}^{1}}$} &
      \vline\nulltab{
        \cball{\csix}{} & \cball{\csix}{} \\
        \cball{\cfiv}{} & & \cball{\csix}{} \\
        \cball{\cfou}{r} & \cball{\cfou}{} \\
        & & \cball{\cfou}{s} & \cball{\cfou}{} & \cball{\cfou}{} & \cball{\cfou}{t} \\
        & \cball{\cfiv}{} & \cball{\cfiv}{} & \cball{\cfiv}{} & \cball{\cfiv}{} & & x \\
        & & & \cball{\csix}{} & & \cball{\cfiv}{} & \cball{\cfou}{} \\\hline }
    \end{array}
  \end{displaymath}
  \caption{\label{fig:exchange}Exchanging labels of an $r$-exchange sequence.}
\end{figure}

It will be convenient to iterate label exchanges as follows; see Fig.~\ref{fig:exchange}.

\begin{definition}
  A sequence of labels $r_0,r_1,\ldots,r_m$ in columns $c_0<c_1<\cdots<c_m$ of a properly labeled diagram is an \newword{$r_0$-exchange sequence} if each $r_i$ abuts $r_{i-1}$ in column $c_{i-1}$ and $r_j$ does not touch $r_i$ in any column $c_{i}<c\le c_j$.
\label{def:exchange-seq}
\end{definition}

\begin{theorem}
  For $r_0,r_1,\ldots,r_m$ in columns $c_0<c_1<\cdots<c_m$ an $r_0$-exchange sequence of a semi-properly labeled diagram, the composition of exchange labelings $\Label_{r_0,r_1}^{c_0} \cdots \Label_{r_{m-1},r_m}^{c_{m-1}}$ is semi-proper, and the $r_i$ in column $c_i$ is re-labeled to $r_0$.
  \label{thm:exchange}
\end{theorem}

\begin{proof}
  We proceed by induction on $m$. The base case $m=0$ is trivial. Assume the theorem for exchange sequences of length less than $m$. Since $r_m$ abuts $r_{m-1}$ at $c_{m-1}$, the exchange labeling $\Label_{r_{m-1},r_m}^{c_{m-1}}$ applies and, by Lemma~\ref{lem:exchange}, results in a semi-proper labeling. Since $r_m$ does not touch $r_{m-1}$ at any column $c\le c_m$, the column $d$ in Definition~\ref{def:exchange-lab} must satisfy $d \ge c_m$, and so the label $r_m$ in column $c_m$ changes to $r_{m-1}$. Since $r_m$ does not touch $r_i$ for $i<m$ in any column $c_i < c \le c_m$, the re-labeled $r_{m-1}$'s still do not touch any $r_i$ for $i<m-1$ in any column $c_i < c \le c_m$. Thus omitting the label $r_m$ and the column $c_{m-1}$ results in an $r_0$-exchange sequence of length $m-1$. The result now follows by induction.  
\end{proof}

\subsection{Labeled insertion}
\label{sec:insert}

\begin{figure}[ht]
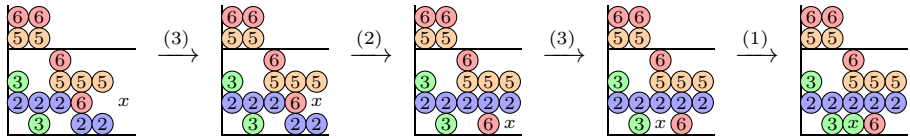

  \begin{displaymath}
    \arraycolsep=4pt
    \begin{array}{ccccccccccc}
      \vline\nulltab{ \cball{\csix}{6} & \cball{\csix}{6} \\
        \cball{\cfiv}{5} & \cball{\cfiv}{5} \\\hline & & \cball{\csix}{6} \\
        \cball{\cthr}{3} & & \cball{\cfiv}{5} & \cball{\cfiv}{5} & \cball{\cfiv}{5} \\
        \cball{\ctwo}{2} & \cball{\ctwo}{2} & \cball{\ctwo}{2} & \cball{\csix}{6} & & x\\
        & \cball{\cthr}{3} & & \cball{\ctwo}{2} & \cball{\ctwo}{2} \\\hline}
      & \raisebox{-1.6\cellsize}{$\xrightarrow{(3)}$} & 
      \vline\nulltab{ \cball{\csix}{6} & \cball{\csix}{6} \\
        \cball{\cfiv}{5} & \cball{\cfiv}{5} \\\hline & & \cball{\csix}{6} \\
        \cball{\cthr}{3} & & \cball{\cfiv}{5} & \cball{\cfiv}{5} & \cball{\cfiv}{5} \\
        \cball{\ctwo}{2} & \cball{\ctwo}{2} & \cball{\ctwo}{2} & \cball{\csix}{6} & x\\
        & \cball{\cthr}{3} & & \cball{\ctwo}{2} & \cball{\ctwo}{2} \\\hline}
      & \raisebox{-1.6\cellsize}{$\xrightarrow{(2)}$} & 
      \vline\nulltab{ \cball{\csix}{6} & \cball{\csix}{6} \\
        \cball{\cfiv}{5} & \cball{\cfiv}{5} \\\hline & & \cball{\csix}{6} \\
        \cball{\cthr}{3} & & \cball{\cfiv}{5} & \cball{\cfiv}{5} & \cball{\cfiv}{5} \\
        \cball{\ctwo}{2} & \cball{\ctwo}{2} & \cball{\ctwo}{2} & \cball{\ctwo}{2} & \cball{\ctwo}{2}\\
        & \cball{\cthr}{3} & & \cball{\csix}{6} & x \\\hline}
      & \raisebox{-1.6\cellsize}{$\xrightarrow{(3)}$} & 
      \vline\nulltab{ \cball{\csix}{6} & \cball{\csix}{6} \\
        \cball{\cfiv}{5} & \cball{\cfiv}{5} \\\hline & & \cball{\csix}{6} \\
        \cball{\cthr}{3} & & \cball{\cfiv}{5} & \cball{\cfiv}{5} & \cball{\cfiv}{5} \\
        \cball{\ctwo}{2} & \cball{\ctwo}{2} & \cball{\ctwo}{2} & \cball{\ctwo}{2} & \cball{\ctwo}{2}\\
        & \cball{\cthr}{3} & x & \cball{\csix}{6} & \\\hline}
      & \raisebox{-1.6\cellsize}{$\xrightarrow{(1)}$} & 
      \vline\nulltab{ \cball{\csix}{6} & \cball{\csix}{6} \\
        \cball{\cfiv}{5} & \cball{\cfiv}{5} \\\hline & & \cball{\csix}{6} \\
        \cball{\cthr}{3} & & \cball{\cfiv}{5} & \cball{\cfiv}{5} & \cball{\cfiv}{5} \\
        \cball{\ctwo}{2} & \cball{\ctwo}{2} & \cball{\ctwo}{2} & \cball{\ctwo}{2} & \cball{\ctwo}{2}\\
        & \cball{\cthr}{3} & \cball{\cthr}{x} & \cball{\csix}{6} & \\\hline}
    \end{array}
  \end{displaymath}
  \caption{\label{fig:rest-rect}Rectification of a labeled diagram with $x$ restricted by $4$.}
\end{figure}

Our insertion procedure is similar to RSK in that an inserted cells may \emph{bump} other cells within its column, though when no bump is possible, the cell may also slide left to pass other cells within its row; see Fig.~\ref{fig:rest-rect}.

Extending prior notation, a position $x\not\in T$ in column $c+1$ \newword{abuts} $r$ if $r$ in column $c$ lies weakly above $x$ and either $\alpha_r=c$ or $r$ in column $c+1$ lies strictly below $x$. For example, in Fig.~\ref{fig:exchange}, $x$ abuts $t$, then $s$, then $r$. 

\begin{definition}
  For $T\in\PKD(\alpha)$ and $x\not\in T$ a vacant position weakly below row $k$, the \newword{rectification of $T$ labeled by $\alpha$ with $x$ restricted by $k$}, denoted by $\Rect^{(k)}_{\alpha}(T,x)$, is defined by 
  \begin{enumerate}
  \item if $T$ has a semi-proper labeling of shape $\alpha'\lswappin\alpha$ in which $x$ abuts the rightmost $r \le k$, then return $T \cup \{x\}$;
  \item else if $T$ has a semi-proper labeling of shape $\alpha'\lswappin\alpha$ in which $x$ abuts $r \le k$,  then return $\Rect_{\alpha}^{(k)}(T\cup\{x\}\setminus\{y\},y)$ for $y$ the lowest cell in the column of $x$ for which some semi-proper labeling $\Label$ has $x$ abut $r \le k$ and $\Label(y) = r$;
  \item else return $\Rect_{\alpha}^{(k)}(T,y)$ for $y$ the rightmost vacant position left of $x$.
  \end{enumerate}
    \label{def:rest-rect}
\end{definition}

We refer to case (2) as a \newword{bump} and to case (3) as a \newword{pass}. The \newword{bumping path} is the set of distinguished cells, and the \newword{landing cell} is the final distinguished cell. 

\begin{theorem}
  For $T\in\PKD(\beta)$ and $x\not\in T$ a vacant position weakly below row $k$, $\Rect^{(k)}_{\beta}(T,x)$ is a well-defined, rectified diagram.    
  \label{thm:well-def}
\end{theorem}

\begin{proof}
  Suppose case (3) occurs, say with $x$ in some row $r$. Since $x$ lies weakly below row $k$, we must have $r \le k$. Suppose, for contradiction, there is no vacant position left of $x$ within its row. From left to right, let $r,r_1,\ldots,r_m$ be the distinct labels, say with each occurring last in column $c_0<c_1<\cdots<c_m$. We claim this is an $r$-exchange sequence, contradicting that $R=\varnothing$. For $i>0$, each $r_i$ abuts $r_{i-1}$ at $c_{i-1}$ since, by Definition~\ref{def:KT}\eqref{i:descend}, any $r_{i-1}$ in column $c_{i-1}+1$ lies weakly lower, and hence lower since $r_i$ is the label in row $r$, column $c_{i-1}+1$. Moreover, for $i<j$, all labels $r_i$ in columns $c_i < c$ lie strictly below row $r$, and all labels $r_j$ in columns $c \le c_j$ lie weakly above row $r$. Thus for $i<j$, $r_j$ cannot touch $r_i$ in any column $c_i < c \le c_j$. We have a contradiction, and so a vacant position exists.

  Suppose case (2) occurs, with $y$ the lowest cell and $r$ its label after an $r$-exchange labeling. Since $x$ abuts $r$, $y$ lies below $x$ and hence below row $k$. Label $x$ with $r$ and, after removing $y$, we claim this is semi-proper. The exchange labeling is semi-proper by Theorem~\ref{thm:exchange}, and since $x$ abuts $r$, it lies weakly below the last $r$ left of $y$. Since $x$ lies above $y$, all subsequent $r$'s lie weakly below it, and so the labeling is semi-proper. By Lemma~\ref{lem:matching}, $T\cup\{x\}\setminus\{y\}$ is a Kohnert diagram for the exchanged labeling $\alpha \lswappin \beta$, and so by Lemma~\ref{lem:lswappin}, $T\cup\{x\}\setminus\{y\}\in\PKD(\beta)$. 

  Suppose case (1) occurs for some $r\le k$. We claim labeling $x$ with $r$ is semi-proper. Indeed, the exchange labeling for the $r$-exchange sequence is semi-proper by Theorem~\ref{thm:exchange}, and since $x$ abuts $r$, it lies weakly below the last $r$. Thus by Lemma~\ref{lem:matching}, $T\cup\{x\}$ is a Kohnert diagram for $\alpha+\e_r$ for some $\alpha \lswappin \beta$. Finally, since the bumping path during rectification moves either strictly left or strictly down, the process must eventually terminate with a rectified diagram.
\end{proof}

\begin{definition}
  For $T\in\KD(\beta)$ and $1 \le r \le k$ integers, the \newword{insertion of $r$ restricted by $k$ into $T$ labeled by $\beta$}, denoted by $T \xleftarrow{\beta,k} r$, is
  \begin{equation}
    T \xleftarrow{\beta,k} r = \Rect^{(k)}_{\alpha}(T,(c+1,r)) .
  \end{equation}
  where $c=\max(\beta_i)$ and $\alpha$ is maximal such that $\thread(T)\lswappin\alpha\lswap\beta$.
  \label{def:rest-insert}
\end{definition}

\begin{figure}[ht]
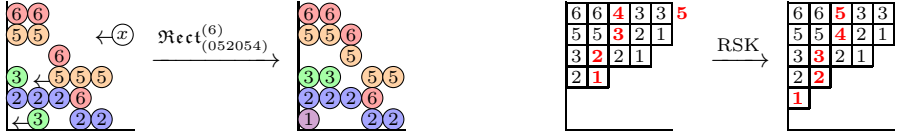

  \begin{displaymath}
    \arraycolsep=0.4\cellsize
    \begin{array}{ccccccc}
      \vline\nulltab{\cball{\csix}{6} & \cball{\csix}{6} \\
        \cball{\cfiv}{5} & \cball{\cfiv}{5} & & & & \leftball{white}{x} \\
        & & \cball{\csix}{6} \\
        \cball{\cthr}{3} & & \leftball{\cfiv}{5} & \cball{\cfiv}{5} & \cball{\cfiv}{5} \\
        \cball{\ctwo}{2} & \cball{\ctwo}{2} & \cball{\ctwo}{2} & \cball{\csix}{6} \\
        & \leftball{\cthr}{3} & & \cball{\ctwo}{2} & \cball{\ctwo}{2}  \\\hline} 
      & \raisebox{-2\cellsize}{$\xrightarrow{\Rect_{(052054)}^{(6)}}$} &
      \vline\nulltab{\cball{\csix}{6} & \cball{\csix}{6} \\
        \cball{\cfiv}{5} & \cball{\cfiv}{5} & \cball{\csix}{6} \\
        & & \cball{\cfiv}{5} \\
        \cball{\cthr}{3} & \cball{\cthr}{3} & & \cball{\cfiv}{5} & \cball{\cfiv}{5} \\
        \cball{\ctwo}{2} & \cball{\ctwo}{2} & \cball{\ctwo}{2} & \cball{\csix}{6} \\
        \cball{\cone}{1} & & & \cball{\ctwo}{2} & \cball{\ctwo}{2} \\\hline}
      & \hspace{6\cellsize} &
      \vline\tableau{
      6 & 6 & \mathbf{\color{red}4} & 3 & 3\\
      5 & 5 & \mathbf{\color{red}3} & 2 & 1 \\
      3 & \mathbf{\color{red}2} & 2 & 1 \\
      2 & \mathbf{\color{red}1} \\
      & \\ &  \\\hline}\hss\nulltab{\mathbf{\color{red}5} \\ \\ \\ \\ \\ }
      & \raisebox{-2\cellsize}{$\xrightarrow{\mathrm{RSK}}$} &
      \vline\tableau{
      6 & 6 & \mathbf{\color{red}5} & 3 & 3 \\
      5 & 5 & \mathbf{\color{red}4} & 2 & 1 \\
      3 & \mathbf{\color{red}3} & 2 & 1 \\
      2 & \mathbf{\color{red}2} \\
      \mathbf{\color{red}1} & \\ & \\\hline} 
    \end{array}
  \end{displaymath}
  \caption{\label{fig:rectification}The insertion of $5$ restricted by $6$ into a diagram $T\in\KD(052054)$ and RSK insertion of $5$ into $\varphi(T)\in\SSRT(002455)$.}
\end{figure}

To persuade the skeptic, note restricted rectification generalizes RSK.; see Fig.~\ref{fig:rectification}.

\begin{theorem}
  For $T\in\KD(\alpha)$ and $k \ge \ell(\alpha)$, we have
  \begin{equation}
    \varphi(T \xleftarrow{\alpha,k} r) = \varphi(T) \xleftarrow{\mathrm{RSK}} r.
  \end{equation}
  where $\varphi$ is the injective, weight-preserving map $\KD(\alpha) \rightarrow \SSRT(\sort(\alpha))$.
  \label{thm:RSK}
\end{theorem}

\begin{proof}
  If $T\cup\{x\}$ is a rectified diagram, then $x$ must abut some label since $x$ can be labeled consistently with $T$. Thus when $k \ge \ell(\alpha)$, case (3) of Definition~\ref{def:rest-rect} occurs if and only if $T\cup\{x\}$ is not a rectified diagram, and in this case there is no cell immediately to the left of $x$. Therefore Definition~\ref{def:rest-rect} reduces to Assaf's original \emph{rectification} \cite[Def~4.2.4]{Ass22-KC}, a procedure that maps an arbitrary diagram to a rectified one of the same weight. The theorem now follows from \cite[Thm~4.2.7]{AQ}, in which Assaf and Quijada prove rectification, in the unrestricted sense, of a diagram with one additional cell appended is equivalent to RSK insertion.
\end{proof}

\begin{figure}[ht]
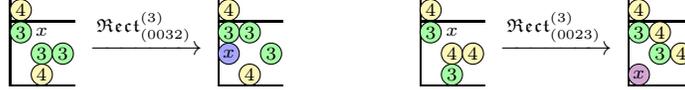

  \begin{displaymath}
    \arraycolsep=3pt
    \begin{array}{ccccccc}
      \vline\nulltab{
        \cball{\cfou}{4} \\ \hline
        \cball{\cthr}{3} & x \\
        & \cball{\cthr}{3} & \cball{\cthr}{3}\\
        & \cball{\cfou}{4} \\\hline} 
      & \raisebox{-1.6\cellsize}{$\xrightarrow{\Rect_{(0032)}^{(3)}}$} & 
      \vline\nulltab{
        \cball{\cfou}{4} \\ \hline
        \cball{\cthr}{3} & \cball{\cthr}{3} \\
        \cball{\ctwo}{x} & & \cball{\cthr}{3}\\
        & \cball{\cfou}{4} \\\hline} 
      & \hspace{5\cellsize} &
      \vline\nulltab{
        \cball{\cfou}{4} \\ \hline
        \cball{\cthr}{3} & x \\
        & \cball{\cfou}{4} & \cball{\cfou}{4}\\
        & \cball{\cthr}{3} \\\hline} 
      & \raisebox{-1.6\cellsize}{$\xrightarrow{\Rect_{(0023)}^{(3)}}$} & 
      \vline\nulltab{
        \cball{\cfou}{4} \\ \hline
        \cball{\cthr}{3} & \cball{\cfou}{4}  \\
        & \cball{\cthr}{3} & \cball{\cfou}{4}\\
        \cball{\cone}{x} & \\\hline} 
    \end{array}
  \end{displaymath}
  \caption{\label{fig:label}An example showing labels matter for rectification.}
\end{figure}

In general, the label used for restricted rectification affects the result; see Fig~\ref{fig:label}.

\subsection{Single cell excision}
\label{sec:inverse}

We can uniquely identify the landing column of $U=(T \xleftarrow{\beta,k} t)$ as the unique column of $U$ with more cells than $\D(\beta)$. To reverse the insertion, we must locate the landing cell within this column. To begin, we observe passing is a last resort.

\begin{lemma}
  Let $U\in\KD(\beta)$, $z\not\in U$ a vacant position (or $z = \varnothing$ regarded in row $0$), and $x\in U$ a cell for which $T = U \cup\{z\} \setminus\{x\} \in\KD(\beta)$ and, for $\alpha$ as in Definition~\ref{def:rest-insert}, suppose $\Rect^{(k)}_{\alpha}(T,x)$ bumps $z$. Let $y$ be any cell of $U$ weakly above $z$ and strictly below $x$ for which $S = U \cup\{z\}\setminus\{y\}$ is rectified. Then for $u\not\in U$ the nearest vacant position to the right of $y$ within its row, there exists $r\le k$ for which $u$ abuts $r$ after exchanging labels for some $r$-exchange sequence. 
  \label{lem:donotpass}
\end{lemma}

\begin{figure}[ht]
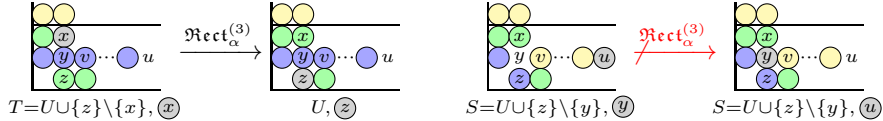

  \begin{displaymath}
    \arraycolsep=-1pt
    \begin{array}{ccccccc}
      \vline\nulltab{
        \cball{\cfou}{} & \cball{\cfou}{} \\\hline
        \cball{\cthr}{} & \cball{gray}{x} \\
        \cball{\ctwo}{} & \cball{\ctwo}{y} & \cball{\ctwo}{v} & \cdots & \cball{\ctwo}{} & u\\
        & \cball{\cthr}{z} & \cball{\cthr}{} \\\hline} 
      & \raisebox{-1.6\cellsize}{$\xrightarrow{\Rect_{\alpha}^{(3)}}$} & 
      \vline\nulltab{
        \cball{\cfou}{} & \cball{\cfou}{} \\\hline
        \cball{\cthr}{} & \cball{\cthr}{x} \\
        \cball{\ctwo}{} & \cball{\ctwo}{y} & \cball{\ctwo}{v} & \cdots & \cball{\ctwo}{} & u\\
        & \cball{gray}{z} & \cball{\cthr}{} \\\hline}
      & \hspace{3\cellsize} &
      \vline\nulltab{
        \cball{\cfou}{} & \cball{\cfou}{} \\\hline
        \cball{\cthr}{} & \cball{\cthr}{x} \\
        \cball{\ctwo}{} & y & \cball{\cfou}{v} & \cdots & \cball{\cfou}{} & \cball{gray}{u}\\
        & \cball{\ctwo}{z} & \cball{\cthr}{} \\\hline}
      & \raisebox{-1.6\cellsize}{$\color{red}\not\xrightarrow{\Rect_{\alpha}^{(3)}}$} & 
      \vline\nulltab{
        \cball{\cfou}{} & \cball{\cfou}{} \\\hline
        \cball{\cthr}{} & \cball{\cthr}{x} \\
        \cball{\ctwo}{} & \cball{gray}{y} & \cball{\cfou}{v} & \cdots & \cball{\cfou}{} & u\\
        & \cball{\ctwo}{z} & \cball{\cthr}{} \\\hline} \\
      \scriptstyle T = U \cup\{z\}\setminus\{x\}, \raisebox{-0.4\cellsize}{$\cball{gray}{x}$} & &
      \hspace{2em} \scriptstyle U, \raisebox{-0.4\cellsize}{$\cball{gray}{z}$} \hspace{2em} & &
      \scriptstyle S = U \cup\{z\}\setminus\{y\}, \raisebox{-0.4\cellsize}{$\cball{gray}{y}$} & & 
      \scriptstyle S = U \cup\{z\}\setminus\{y\}, \raisebox{-0.4\cellsize}{$\cball{gray}{u}$}
    \end{array}
  \end{displaymath}
  \caption{\label{fig:donotpass}An illustration for the No Passing Lemma~\ref{lem:donotpass}.}
\end{figure}

\begin{proof}
  Consider $U$ with its rectified labeling, and let $r \le k$ be the label assigned to $x$. Since $\Rect^{(k)}_{\alpha}(T,x)$ bumps $z$ (or lands), there must be an $r$ in the column immediately to the right of $z$ and, by Definition~\ref{def:KT}\eqref{i:descend}, it lies weakly below the row of $z$ (or no $r$'s strictly to the right of $x$). Let $v$ be the cell of $S$ immediately to the right of the vacant position $y$, which lies weakly left of $u$; see Fig.~\ref{fig:donotpass}. Since $U\setminus\{y\}$ is rectified and $y$ lies above $z$, we may label $S$ the same as $U$ except for the cells right of $y$ that share its label, which might be re-labeled. In particular, the cells with labels $r$ in $U$ still have label $r$ in $S$, and so $v$ since $v$ lies below $x$ and weakly above $z$, it abuts $r$. Since $u$ lies in the same row as $v$ with no empty positions between, this ensures $u$ abuts $r$ after an $r$-exchange sequence, as desired. 
\end{proof}

Using Lemma~\ref{lem:donotpass}, we now prove restricted insertion is reversible.

\begin{theorem}
  If $(S \xleftarrow{\beta,k} s) = (T \xleftarrow{\beta,k} t)$ for some $S,T\in\KD(\beta)$ and some $1 \le s,t\le k$, then $S=T$ and $s=t$.
  \label{thm:inverse}
\end{theorem}

\begin{proof}
  By Proposition~\ref{prop:bump}, both insertions have the same landing column. If either insertion results in an immediate landing (case (1)), then the landing cell lies in column $\max(\beta_i)+1$, and so the insertion rows must coincide. Thus $s=t$. Since the final diagrams are the same and no adjustments were made to $S$ or to $T$, this also forces $S=T$ as well. Since the original insertion occurs in a column not occupied by any cell of $S,T$, the first nontrivial step in rectification must be a pass. In particular, both insertions include a pass at some point.

  Suppose rectifications for $S$ and $T$ have different landing cells within the same column. Since both rectifications must include a pass at some point, we may consider the cells $x_S$ and $x_T$ which mark the first position within the landing column for $S$ and $T$, respectively. If $x_S = x_T$, then since $(S \xleftarrow{\beta,k} s) = (T \xleftarrow{\beta,k} t)$ and this is the landing column, both would have the same landing cell. Thus none of the bumped cells within this column may coincide. Without loss of generality, say $x_S$ lies below $x_T$. We may continue with the rectification of $T$ until the point at which a cell lands above $x_S$ or bumps a cell below $x_S$. By Lemma~\ref{lem:donotpass}, this contradicts the possibility that $x_S$ entered the column by a pass. Thus both insertions have the same landing cell, and so $S,T\in\PKD(\alpha)$ for the same rectification label $\alpha$.

  Reading the rectification process in reverse, consider the diagram $U$ with vacant position $x_U$ that marks the last position before the two processes diverge. That is, assume for diagrams $S,T\in\PKD(\alpha)$, both properly labeled of shape $\alpha$, and vacant positions $x_S\not\in S$ and $x_T\not\in T$, either $x_S \neq x_T$ or $S \neq T$ (or both) and one nontrivial step (a bump or a pass) of the rectification process for $\Rect^{(k)}_{\alpha}(S,x_S)$ and for $\Rect^{(k)}_{\alpha}(T,x_T)$ results in the diagram $U$ with vacant position $x_U$. 

  If both $\Rect^{(k)}_{\alpha}(S,x_S)$ and $\Rect^{(k)}_{\alpha}(T,x_T)$ result in a pass, then both $x_S,x_T$ must be the nearest vacant position not in $U$ to the right of $x_U$ within its row, so $x_S=x_T$ and $S = U \cup \{x_S\} \setminus \{x_U\} = U \cup \{x_T\} \setminus \{x_U\} = T$, a contradiction to this being the point of divergence. Thus at least one of the steps must have been a bump.

  By Lemma~\ref{lem:donotpass}, we cannot have one of the two passing to $x_U$ while the other bumps, so both must bump, and so both $x_S,x_T$ lie above $x_U$ within its column. If $x_T = x_S$, then $S=T$ as well, so this isn't the case now, or in any preceding step that also resulted in a bump. Since both rectifications must include a pass at some point, we may consider the lowest cell, say coming from $S$, that entered the column by a pass. Rename this cell $x_S$, and name $x_T$ the nearest position above this which came from $T$. Then $x_T$ must bump a cell strictly below $x_S$, so Lemma~\ref{lem:donotpass} once again applies, contradicting that $x_S$ entered the column from a pass. All cases result in contradiction, so we must have $x_S=x_T$ and $S=T$ throughout.
\end{proof}

\vspace{-0.5\baselineskip}

\begin{corollary}
  If $U = T \xleftarrow{\beta,k} r$ for some $T\in\KD(\beta)$ and some $r \le k$, then, for $\alpha$ as in Def~\ref{def:rest-insert}, the landing cell $x$ is the highest cell weakly below row $k$ in the landing column for which $U \setminus \{x\}\in\KD(\alpha)$ and $\Rect^{(k)}_{\alpha}(\Label_{\alpha}(U \setminus \{x\}),x)$ lands.
  \label{cor:land}
\end{corollary}

In summary, given $\beta$, $k$ and the (unlabeled) diagram $U = T \xleftarrow{\beta,k} r$, we uniquely recover $T$ and $r$ as follows: the landing column is the unique part of $\thread(U)$ that occurs with smaller multiplicity in $\beta$, the landing cell $x$ is the highest cell weakly below row $k$ for which $\Rect^{(k)}_{\alpha}(\Label_{\alpha}(U \setminus \{x\}),x)$ lands, and we must pass $x$ to the right to nearest vacant position not in $U$. Continuing, replace $x$ with the highest cell within its column which could have bumped $x$, or, if none is found, pass $x$ the right. The process terminates once the vacant position passes into a column with no cells, and this is the row $r$, now excised from the diagram $T$.


%
\section{Iterated insertion}
%
\label{sec:atomexpansion}

In this section, we iterate insertions into a Kohnert diagram, which leads to saturated chains of compositions of increasing size.

\subsection{Re-labeling after insertion}
\label{sec:iterate}

To multiply by a Schur polynomial, we must consider successive insertions, each with the same parameter $k$, determined by the Schur polynomial, though the labels must change as the diagram grows.

\begin{definition}
  Write \newword{$\beta\kup_{k}t_{r,s_m} \cdots t_{r,s_1} \beta + \e_r$} whenever $r \le k < s_1 < \ldots < s_m$ and $0 < \beta_{r} < \beta_{s_1} < \cdots <\beta_{s_m}$. Denote the transitive closure by \newword{$\beta\kup_{k}\gamma$}.
  \label{def:kuppin}
\end{definition}

\begin{figure}[ht]
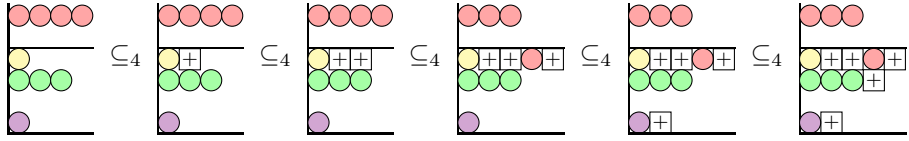

  \begin{displaymath}
    \arraycolsep=3pt
    \begin{array}{ccccccccccc}
      \vline\nulltab{
        \cball{\csix}{} & \cball{\csix}{} & \cball{\csix}{} & \cball{\csix}{} \\ \\\hline
        \cball{\cfou}{} & \\
        \cball{\cthr}{} & \cball{\cthr}{} & \cball{\cthr}{} \\
        & \\
        \cball{\cone}{} \\\hline}
      & \raisebox{-1.7\cellsize}{$\kup_4$} &
      \vline\nulltab{
        \cball{\csix}{} & \cball{\csix}{} & \cball{\csix}{} & \cball{\csix}{} \\ \\\hline
        \cball{\cfou}{} & \cbox{white}{+} \\
        \cball{\cthr}{} & \cball{\cthr}{} & \cball{\cthr}{} \\
        & \\
        \cball{\cone}{} \\\hline}
      & \raisebox{-1.7\cellsize}{$\kup_4$} &
      \vline\nulltab{
        \cball{\csix}{} & \cball{\csix}{} & \cball{\csix}{} & \cball{\csix}{} \\ \\\hline
        \cball{\cfou}{} & \cbox{white}{+} & \cbox{white}{+} \\
        \cball{\cthr}{} & \cball{\cthr}{} & \cball{\cthr}{}  \\
        & \\
        \cball{\cone}{} \\\hline}
      & \raisebox{-1.7\cellsize}{$\kup_4$} &
      \vline\nulltab{
        \cball{\csix}{} & \cball{\csix}{} & \cball{\csix}{}  \\ \\\hline
        \cball{\cfou}{} & \cbox{white}{+} & \cbox{white}{+} & \cball{\csix}{} & \cbox{white}{+} \\
        \cball{\cthr}{} & \cball{\cthr}{} & \cball{\cthr}{} \\
        & \\
        \cball{\cone}{} \\\hline}
      & \raisebox{-1.7\cellsize}{$\kup_4$} &
      \vline\nulltab{
        \cball{\csix}{} & \cball{\csix}{} & \cball{\csix}{}  \\ \\\hline
        \cball{\cfou}{} & \cbox{white}{+} & \cbox{white}{+} & \cball{\csix}{} & \cbox{white}{+} \\
        \cball{\cthr}{} & \cball{\cthr}{} & \cball{\cthr}{} \\
        & \\
        \cball{\cone}{} & \cbox{white}{+} \\\hline}
      & \raisebox{-1.7\cellsize}{$\kup_4$} &
      \vline\nulltab{
        \cball{\csix}{} & \cball{\csix}{} & \cball{\csix}{}  \\ \\\hline
        \cball{\cfou}{} & \cbox{white}{+} & \cbox{white}{+} & \cball{\csix}{} & \cbox{white}{+} \\
        \cball{\cthr}{} & \cball{\cthr}{} & \cball{\cthr}{} & \cbox{white}{+} \\
        & \\
        \cball{\cone}{} & \cbox{white}{+} \\\hline}
    \end{array} 
  \end{displaymath}
    \caption{\label{fig:chain}A saturated chain in $\kup_4$ from $\beta=(1,0,3,1,0,4)$ to $\gamma=(2,0,4,5,0,3)$ with added columns $2,3,5,2,4$ and extended rows $4,4,4,1,3$, respectively.}
\end{figure}

That is, we apply left swaps with $r \le k < s_i$ and then increment the length of row $r$. We call $r$ the \newword{extended row} and $\beta_{s_m}+1$ the \newword{added column}; see Fig.~\ref{fig:chain}.

\begin{proposition}
  For $T\in\PKD(\beta)$ and $x\not\in T$ a vacant position weakly below $k$, we have $\Rect_{\beta}^{(k)}(T,x)\in\PKD(\gamma)$ for some $\beta\kupdot_{k}\gamma$.
  \label{prop:bump}
\end{proposition}

\begin{proof}
  Let $U = \Rect_{\beta}^{(k)}(T,x)$, and let $y$ denote the landing cell and $c$ the landing column. In any label exchange during rectification, the first column is unaffected and every affected column has the larger label. Thus the change in the shape of the labeling, if any, is a pinned left swap, and so, by Lemma~\ref{lem:lswappin}, $U\setminus\{y\}\in\PKD(\beta)$.

  Since $y$ lands, it abuts some label $r\le k$ for which there is no label $r$ in column $c$ even after applying any necessary $r$-exchange sequence. We claim there is also no $r_i$ in column $c$ for any label $r_i$ on the $r$-exchange sequence. Indeed, if $r_i$ is the largest label on the $r$-exchange sequence for which there is an $r_i$ in column $c$, then this $r_i$ must be changed to $r_{i-1}$, and so on, until at last there is an $r$ in column $c$, contradicting that this is the landing column. Thus the claim is proved.

  There is an $r$-exchange sequence with all other labels greater than $k$, since if $r < s \le k$, then by the previous claim there is no $s$ in column $c$, so $y$ abuts $s$ via the $s$-exchange subsequence. Thus we have an $r$-exchange sequence with $r \le k < r_1 < \cdots < r_m$ where $\beta_r < \beta_{r_1} \le \cdots \le \beta_{r_m}$. Omitting transposition which act trivially, after exchanging labels the shape becomes
  \[ t_{r,r_1} t_{r_1,r_2} \cdots t_{r_{m-1},r_m} \beta = t_{r,r_m} \cdots t_{r,r_2}  t_{r,r_1} \beta \lswappin_{k} \beta . \]
  Allowing $y$ to take label $r$ adds $\e_r$, and the result follows from Lemma~\ref{lem:matching}.
\end{proof}

\begin{figure}[ht]
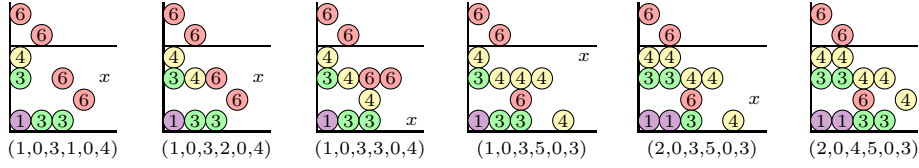

  \begin{displaymath}
    \arraycolsep=\cellsize
    \begin{array}{cccccc}
      \vline\nulltab{
        \cball{\csix}{6} \\ & \cball{\csix}{6} \\\hline
        \cball{\cfou}{4} \\
        \cball{\cthr}{3} & & \cball{\csix}{6} & & x\\
        & & & \cball{\csix}{6}  \\
        \cball{\cone}{1} & \cball{\cthr}{3} & \cball{\cthr}{3} &  \\\hline}
      &
      \vline\nulltab{
        \cball{\csix}{6} \\ & \cball{\csix}{6} \\\hline
        \cball{\cfou}{4} \\
        \cball{\cthr}{3} & \cball{\cfou}{4} & \cball{\csix}{6} & & x \\
        & & & \cball{\csix}{6} & \\
        \cball{\cone}{1} & \cball{\cthr}{3} & \cball{\cthr}{3} &  \\\hline}
      &
      \vline\nulltab{
        \cball{\csix}{6} \\ & \cball{\csix}{6} \\\hline
        \cball{\cfou}{4} & & \\
        \cball{\cthr}{3} & \cball{\cfou}{4} & \cball{\csix}{6} & \cball{\csix}{6} \\
        & & \cball{\cfou}{4} \\
        \cball{\cone}{1} & \cball{\cthr}{3} & \cball{\cthr}{3} & & x \\\hline}
      &
      \vline\nulltab{
        \cball{\csix}{6} \\ & \cball{\csix}{6} \\\hline
        \cball{\cfou}{4} & & & & & x\\
        \cball{\cthr}{3} & \cball{\cfou}{4} & \cball{\cfou}{4} & \cball{\cfou}{4} \\
        & & \cball{\csix}{6} \\
        \cball{\cone}{1} & \cball{\cthr}{3} & \cball{\cthr}{3} & & \cball{\cfou}{4}  \\\hline}
      &
      \vline\nulltab{
        \cball{\csix}{6} \\ & \cball{\csix}{6} \\\hline
        \cball{\cfou}{4} & \cball{\cfou}{4} & \\
        \cball{\cthr}{3} & \cball{\cthr}{3} & \cball{\cfou}{4} & \cball{\cfou}{4} \\
        & & \cball{\csix}{6} & & & x\\
        \cball{\cone}{1} & \cball{\cone}{1} & \cball{\cthr}{3} & & \cball{\cfou}{4}  \\\hline}
      &
      \vline\nulltab{
        \cball{\csix}{6} \\ & \cball{\csix}{6} \\\hline
        \cball{\cfou}{4} & \cball{\cfou}{4} & \\
        \cball{\cthr}{3} & \cball{\cthr}{3} & \cball{\cfou}{4} & \cball{\cfou}{4} \\
        & & \cball{\csix}{6} & & \cball{\cfou}{4} \\
        \cball{\cone}{1} & \cball{\cone}{1} & \cball{\cthr}{3} & \cball{\cthr}{3}  \\\hline}
      \\
      \scriptstyle (1,0,3,1,0,4) & 
      \scriptstyle (1,0,3,2,0,4) &  
      \scriptstyle (1,0,3,3,0,4) &  
      \scriptstyle (1,0,3,5,0,3) &  
      \scriptstyle (2,0,3,5,0,3) &  
      \scriptstyle (2,0,4,5,0,3) 
    \end{array}
  \end{displaymath}
  \caption{\label{fig:insertion}Iterated labeled insertion restricted by $4$ for $r=3,3,1,4,2$.}
\end{figure}

Proposition~\ref{prop:bump} guides how we relabel cells after insertion; see Fig.~\ref{fig:insertion}. We have the following definitions from \cite[Def~3.2.6]{AQ} and \cite[Def~3.2.3]{AQ}.

\begin{definition}[\cite{AQ}]
  A position $(c,r)$ with $r\le k$ is \newword{$k$-addable for $\beta$} if
  \begin{itemize}
  \item $\beta_r < c$ and if $\beta_r < c-1$, then there exists some $t>k$ such that $\beta_t = c-1$;
  \item for all $r < s \leq k$, either $\beta_s < \beta_r$ or $\beta_s \geq c$.
  \end{itemize}
  For a $k$-addable position, the \newword{$k$-addition of $(c,r)$ to $\beta$} is
  \[ \beta \kplus (c,r) = t_{r_0,r_q} \cdots t_{r_0,r_1} \beta + \e_r, \]
  where $r = r_0 < r_1 < \cdots < r_q$ is the unique sequence of row indices such that $\beta_{r_{0}} < \cdots < \beta_{r_q}=c-1$ and if $r_{i-1} < s < r_i$, then $\beta_s \leq \beta_{r_{i-1}}$ or $\beta_s > \beta_{r_i}$.
  \label{def:kadd}
\end{definition}

Informally, we $k$-add a cell to row $r \le k$ in column $c$ provided we may apply Kohnert moves dropping cells above row $k$ down to row $r$ until it has length $c-1$.

\begin{figure}[ht]
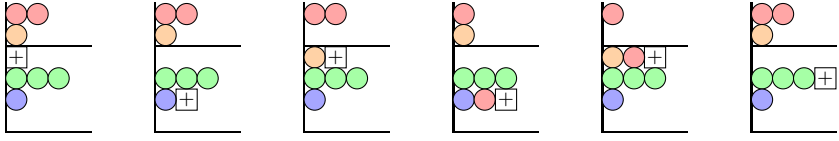

  \begin{displaymath}
    \arraycolsep=1.5\cellsize
    \begin{array}{cccccc}
      \vline\nulltab{
        \cball{\csix}{} & \cball{\csix}{} \\
        \cball{\cfiv}{} \\\hline
        \cbox{white}{+} \\
        \cball{\cthr}{} & \cball{\cthr}{} & \cball{\cthr}{} & \\
        \cball{\ctwo}{} \\
        \\\hline} &
      \vline\nulltab{
        \cball{\csix}{} & \cball{\csix}{} \\
        \cball{\cfiv}{} \\\hline
        \\
        \cball{\cthr}{} & \cball{\cthr}{} & \cball{\cthr}{} & \\
        \cball{\ctwo}{} & \cbox{white}{+} \\
        \\\hline} &
      \vline\nulltab{
        \cball{\csix}{} & \cball{\csix}{} \\
        \\\hline
        \cball{\cfiv}{} & \cbox{white}{+} \\
        \cball{\cthr}{} & \cball{\cthr}{} & \cball{\cthr}{} & \\
        \cball{\ctwo}{} \\
        \\\hline} &
      \vline\nulltab{
        \cball{\csix}{} \\
        \cball{\cfiv}{}  \\\hline
        \\
        \cball{\cthr}{} & \cball{\cthr}{} & \cball{\cthr}{} & \\
        \cball{\ctwo}{} & \cball{\csix}{} & \cbox{white}{+}  \\
        \\\hline} &
      \vline\nulltab{
        \cball{\csix}{} \\
        \\\hline
        \cball{\cfiv}{} & \cball{\csix}{} & \cbox{white}{+}  \\
        \cball{\cthr}{} & \cball{\cthr}{} & \cball{\cthr}{} & \\
        \cball{\ctwo}{} \\
        \\\hline} &
      \vline\nulltab{
        \cball{\csix}{} & \cball{\csix}{} \\
        \cball{\cfiv}{} \\\hline
        \\
        \cball{\cthr}{} & \cball{\cthr}{} & \cball{\cthr}{} & \cbox{white}{+} \\
        \cball{\ctwo}{} \\
        \\\hline} 
    \end{array} 
  \end{displaymath}
  \caption{\label{fig:addable}The six possible $4$-additions to the composition $(013012)$.}
\end{figure}

In particular, notice $\beta \kupdot_{k} \beta\kplus(c,r)$ for any $k$-addable position; see Fig.~\ref{fig:addable}.

\begin{lemma}
  For $T\in\KD(\beta)$ and $r\le k$, we have
  \[ (T \xleftarrow{\beta,k} r) \in \KD(\beta \kplus (c,s)) \]
  for $c$ the landing column and $(c,s)$ a $k$-addable position for $\beta$.
  \label{lem:iterate}
\end{lemma}

\begin{proof}
  Let $U=T \xleftarrow{\beta,k} r$. By Proposition~\ref{prop:bump}, $\thread(U) \lswap \alpha + \e_j$ for some $\alpha\lswap\beta$ and $j\le k$. Thus Lemma~\ref{lem:lswap}, $U\in\KD(\alpha+\e_j)$. By \cite[Thm~3.2.10]{AQ}, we have
  \begin{equation}
    \bigcup_{\alpha\lswap\beta, \ j\le k} \KD(\alpha+\e_j) =
    \bigcup_{(c,r) \ k-\text{addable}} \KD(\beta \kplus (c,r)).
    \label{e:monkey-bijection2}
  \end{equation}
  In particular, $\alpha+\e_j\lswap\beta\kplus(c,s)$ for some $k$-addable position $(c,s)$ for $\beta$. By transitivity, we have $\thread(U)\lswap\beta\kplus(c,s)$. Therefore, by Lemma~\ref{lem:lswap}, there is at least one $k$-addable position $(c,s)$ for $\beta$ for which $U\in\KD(\beta\kplus(c,s))$.
\end{proof}

As in Fig.~\ref{fig:addable}, there may be multiple $k$-addable positions for a given column.

\begin{definition}
  For $T\in\KD(\beta)$ and $r_1,\ldots,r_m\le k$, the \newword{iterated insertion of $r_1,\ldots,r_m$ restricted by $k$ into $T$ labeled by $\beta$}, denoted by $T \xleftarrow{\beta,k} r_1,\ldots,r_m$, is
  \begin{equation}
    (\cdots((T \xleftarrow{\beta,k} r_1) \xleftarrow{\beta^{(1)},k} r_2) \cdots )\xleftarrow{\beta^{(m-1)},k} r_m,
    \label{e:iterate}
  \end{equation}
  where the \newword{rectified label} is $\beta^{(i)} = \beta^{(i-1)} \kplus (c_i,s_i)$ for $c_i$ the landing column of $r_i$ and $s_i$ the \emph{minimal} row index for which $\thread(T \xleftarrow{\beta,k} r_1,\ldots,r_i)\lswap\beta^{(i-1)}\kplus(c_i,s_i)$.
  \label{def:iterate}
\end{definition}

\begin{theorem}
  For $T\in\KD(\beta)$ and $r_1,\ldots,r_m\le k$, the iterated insertion
  \[ U = (T \xleftarrow{\beta,k} r_1,\ldots,r_m)  \]
  is a well-defined, rectified diagram and $\alpha\kup_{k}\alpha^{(m)}$ where $\alpha$ and $\alpha^{(m)}$ are maximal such that $\thread(T)\lswappin\alpha\lswap\beta$ and $\thread(U)\lswappin\alpha^{(m)}\lswap\beta^{(m)}$, respectively.
  \label{thm:iterate}
\end{theorem}

\begin{proof}
  Let $U_i$ denote successive insertion of $r_1,\ldots,r_i$. By Lemma~\ref{lem:iterate}, $U_i\in\KD(\beta^{(i)})$ for $i\ge 1$, and so by Theorem~\ref{thm:well-def}, each $U_i$ is a well-defined rectified diagram. By Proposition~\ref{prop:bump}, there exists $\alpha'\lswappin_{k}\alpha$ and $s\le k$ such that $U_1\in\PKD(\alpha'+\e_s)$. Choose $s$ with $\alpha'_s-\alpha_s$ is minimal, and set $\alpha^{(1)}=\alpha'+\e_s$. By definition, we have $\alpha\kupdot_{k}\alpha^{(1)}$. By Lemma~\ref{lem:lswappin}, $\thread(U_1)\lswappin\alpha^{(1)}$, and by \eqref{e:monkey-bijection2}, $\alpha^{(1)}\lswap\beta^{(1)}$. Moreover, the choice of $s$ and $\beta^{(1)}$ coincide, ensuring $\alpha^{(1)}\lswap\beta^{(1)}$ is the maximal composition with these two properties. Therefore, by Definition~\ref{def:rest-insert}, $U_2 = \Rect_{\alpha^{(1)}}^{(k)}(U_1,(c,r_2))$. Thus we may iterate this argument, obtaining a sequence $\alpha^{(i-1)}\kupdot_{k}\alpha^{(i)}$ for which $\thread(U_i)\lswappin\alpha^{(i)}$ and $\alpha^{(i)}\lswap\beta^{(i)}$.
\end{proof}

This necessary condition on the threads of diagrams is also sufficient.

\begin{lemma}
  If $\beta\kup_{k}\gamma$ and $U\in\KD(\gamma)$, then $U = (T\xleftarrow{\beta,k}r_1,\ldots,r_m)$ for some $T\in\KD(\beta)$ and some $r_1,\ldots,r_m\le k$.
  \label{lem:atomic}
\end{lemma}

\begin{proof}
  We proceed by induction on $m = |\gamma|-|\beta|$, noting the case $m=0$ is trivial. Suppose $\beta \kup_{k} \beta' \kupdot_{k} \gamma$. By Lemma~\ref{lem:matching}, since $U\in\KD(\gamma)$, $U$ can be properly labeled by $\gamma$. Let $r\le k$ such that $\gamma - \e_{r} \lswappin \beta'$, and let $x$ be the rightmost cell of $\Label_{\gamma}(U)$ with label $r$. Set $U' = U \setminus\{x\}$. We claim $U'\in\KD(\beta')$. The labeling on $U'$ inherited from $U$ is semi-proper of shape $\gamma-\e_{r}$ since $x$ was the \emph{rightmost} cell with label $r$. Thus by Lemma~\ref{lem:matching}, $U'\in\KD(\gamma-\e_{r})$, and by Lemma~\ref{lem:lswappin}, $\KD(\gamma-\e_{r})\subseteq\KD(\beta')$, proving the claim.

  By Lemma~\ref{lem:matching}, there exists a proper labeling of $U'$ of shape $\beta'$. Moreover, in this labeling $x$ abuts $r$ through the same label exchanges by which $\gamma - \e_{r} \lswappin \beta'$, and no cell weakly to the right of $x$ has label $r$. Thus $\Rect^{(k)}_{\beta'}(U',x) = U$.

  Therefore we may take $y$ to be the highest cell in the column of $x$, weakly below row $k$ for which $U\setminus\{y\}\in\KD(\beta')$ and $\Rect^{(k)}_{\beta'}(U\setminus\{y\},y) = U$. Now pass $y$ to the right to nearest vacant position not in $U\setminus\{y\}$. Continuing, replace $y$ with the highest cell within its column which could have bumped $y$, or, if none is found, pass $y$ to the right. The process terminates once the vacant position passes into a column with no cells, and this is the row $r'\le k$, now excised from a diagram $T'\in\KD(\beta')$ for which, by Corollary~\ref{cor:land}, $U = T'\xleftarrow{\beta',k}r'$. By induction, there exists $T\in\KD(\beta)$ and $r_1,\ldots,r_{m-1}\le k$ such that $T' = (T\xleftarrow{\beta,k}r_1,\ldots,r_{m-1})$. The theorem follows by inserting $r'$ into $T'$ to obtain $U$.
\end{proof}

It remains to gather the objects in the image into Demazure atoms.

\subsection{Row-bumping lemma}
\label{sec:row-bump}

\begin{figure}[ht]
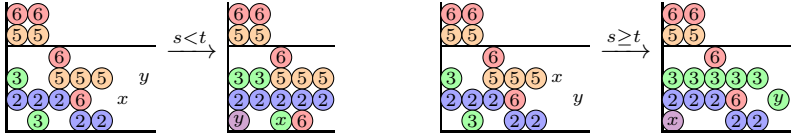

  \begin{displaymath}
    \arraycolsep=2pt
    \begin{array}{ccccccc}
      \vline\nulltab{ \cball{\csix}{6} & \cball{\csix}{6} \\
        \cball{\cfiv}{5}  & \cball{\cfiv}{5} \\\hline &  & \cball{\csix}{6} \\
        \cball{\cthr}{3} & & \cball{\cfiv}{5}  & \cball{\cfiv}{5}  & \cball{\cfiv}{5}  & & y\\
        \cball{\ctwo}{2} & \cball{\ctwo}{2} & \cball{\ctwo}{2} & \cball{\csix}{6} & & x\\
        & \cball{\cthr}{3} & & \cball{\ctwo}{2} & \cball{\ctwo}{2} \\\hline}
      & \raisebox{-1.6\cellsize}{$\xrightarrow{s<t}$} & 
      \vline\nulltab{ \cball{\csix}{6} & \cball{\csix}{6} \\
        \cball{\cfiv}{5}  & \cball{\cfiv}{5} \\\hline &  & \cball{\csix}{6} \\
        \cball{\cthr}{3} & \cball{\cthr}{3} & \cball{\cfiv}{5}  & \cball{\cfiv}{5}  & \cball{\cfiv}{5}  \\
        \cball{\ctwo}{2} & \cball{\ctwo}{2} & \cball{\ctwo}{2} & \cball{\ctwo}{2} & \cball{\ctwo}{2}\\
        \cball{\cone}{y} & & \cball{\cthr}{x} & \cball{\csix}{6} & \\\hline}
      & \hspace{4\cellsize} &
      \vline\nulltab{ \cball{\csix}{6} & \cball{\csix}{6} \\
        \cball{\cfiv}{5}  & \cball{\cfiv}{5} \\\hline &  & \cball{\csix}{6} \\
        \cball{\cthr}{3} & & \cball{\cfiv}{5}  & \cball{\cfiv}{5}  & \cball{\cfiv}{5}  & x\\
        \cball{\ctwo}{2} & \cball{\ctwo}{2} & \cball{\ctwo}{2} & \cball{\csix}{6} & & & y\\
        & \cball{\cthr}{3} & & \cball{\ctwo}{2} & \cball{\ctwo}{2} \\\hline}
      & \raisebox{-1.6\cellsize}{$\xrightarrow{s\ge t}$} & 
      \vline\nulltab{ \cball{\csix}{6} & \cball{\csix}{6} \\
        \cball{\cfiv}{5}  & \cball{\cfiv}{5} \\\hline &  & \cball{\csix}{6} \\
        \cball{\cthr}{3} & \cball{\cthr}{3} & \cball{\cthr}{3} & \cball{\cthr}{3} & \cball{\cthr}{3} \\
        \cball{\ctwo}{2} & \cball{\ctwo}{2} & \cball{\ctwo}{2} & \cball{\csix}{6} & & \cball{\cthr}{y}\\
        \cball{\cone}{x} & & & \cball{\ctwo}{2} & \cball{\ctwo}{2} \\\hline}
    \end{array}
  \end{displaymath}
  \caption{\label{fig:rowbump}Insertion of $2$ then $3$ (left) and $3$ then $2$ (right) into a pinned, properly labeled diagram with restriction parameter $4$.}
\end{figure}

To understand iterated insertions, we have the following analog of Schensted's row bumping lemma \cite{Sch61}; see Fig.~\ref{fig:rowbump}.

\begin{lemma}
  For $S\in\PKD(\alpha)$ and $s,t\ge 1$, let $T = S \xleftarrow{\alpha,k} s$ with landing cell $x$ and rectified labeling $\beta$, and let $U = T \xleftarrow{\beta,k} t$ with landing cell $y$.
  \begin{enumerate}[label=(\roman*)]
  \item\label{i:inc} If $s < t$, then $x$ lies weakly right of $y$ and below $y$ if in the same column.
  \item\label{i:dec} If $s \ge t$, then $x$ lies strictly left of $y$ and weakly above if in adjacent columns.
  \end{enumerate}
  \label{lem:rowbump}
\end{lemma}

\begin{proof}
  Suppose $s < t$, and follow the bumping paths of $x$ for $S$ and $y$ for $T$, which move left and down. If $x$ lands in column $\max(\alpha)+1$, then since $y$ is inserted above, $y$ passes to column $\max(\alpha)+1$, at which point it lands above $x$ or strictly to the left if it bumps. Thus we may assume $x$ passed at some point. Following the path for $y$, suppose, for contradiction, there is some column strictly right of the landing column for $x$ in which either $y$ lands or bumps a cell $z$ below the bumping path for $x$. Take $c$ to be the rightmost such column. Then $y$ must abut some $r\le k$ in column $c$, possibly after an $r$-exchange sequence. Since $x$ lands strictly to the left of $c$, $x$ also encountered this column, and either passed through or bumped out of it. However, since the same $r$ was available for $x$ and the bumped cell $z$ lies below the exit row for $x$, $x$ could have used this $r$-exchange sequence to land or to bump $z$ as well, contradicting that $x$ lands strictly to the left. Therefore the bumping path for $y$ stays above that for $x$ until the landing column of $x$, after which $y$ eventually land in a column weakly to the left of $x$ and above $x$ if in the same column. 

  Suppose $s \ge t$, and again follow the paths. If $x$ lands in column $\max(\alpha)+1$, then $y$ abuts $x$ and so lands in column $\max(\alpha)+2$, which is strictly to the right. Otherwise, suppose, for contradiction, there is some column in which $x$ bumps $z$, say in row $r$, or lands and $y$ passes into or through the column in row $s$ with $s>r$ if $x$ bumps $z$. Take $c$ to be the rightmost such column. Since to the right of $c$ the path for $x$ lies strictly above the path for $y$, the path for $x$ entered column $c$ above the path for $y$. In particular, since $z$, if it exists, lies below the path for $y$, the cell in column $c+1$, row $s$ (or the prior position of $y$ if it passed from column $c+1$) abuts $x$ at column $c$, and so $y$ did not pass. Therefore the bumping path for $y$ stays strictly below and lands in a column strictly to the right of the path for $x$. 
\end{proof}

\subsection{Recoding tableaux}
\label{sec:record}

As in RSK, we use a recording procedure to track the order in which cells were added. For this, we require skew diagrams.

\begin{definition}
  Suppose $\beta\kup_{k}\gamma$ results from the $k$-addition of cells in columns $c_1,\ldots,c_m$ extending rows $r_1,\ldots,r_m$, respectively. For $\alpha\lswap\gamma$, the \newword{skew diagram of shape $\alpha/\beta$} is set of cells $\{x_1,\ldots,x_m\}\subseteq\D(\alpha)$ where $x_i$ is the cell with label $r_i$ in column $c_i$ of $\Label_{\gamma}(\D(\alpha))$. We draw this with non-skewed cells labeled by $\Label_{\beta}$.
  \label{def:skew1}
\end{definition}

\begin{figure}[ht]
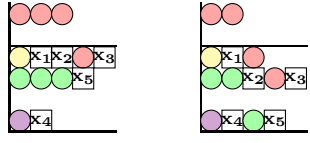

  \begin{displaymath}
    \arraycolsep=2\cellsize
    \begin{array}{cc}
      \vline\nulltab{
        \cball{\csix}{} & \cball{\csix}{} & \cball{\csix}{}  \\ \\\hline
        \cball{\cfou}{} & \cbox{white}{x_1} & \cbox{white}{x_2} & \cball{\csix}{} & \cbox{white}{x_3} \\
        \cball{\cthr}{} & \cball{\cthr}{} & \cball{\cthr}{} & \cbox{white}{x_5} \\
        & \\
        \cball{\cone}{} & \cbox{white}{x_4} \\\hline} &
      \vline\nulltab{
        \cball{\csix}{} & \cball{\csix}{}   \\ \\\hline
        \cball{\cfou}{} & \cbox{white}{x_1} & \cball{\csix}{} \\ 
        \cball{\cthr}{} & \cball{\cthr}{} & \cbox{white}{x_2} & \cball{\csix}{} & \cbox{white}{x_3} \\
        & \\
        \cball{\cone}{} & \cbox{white}{x_4} & \cball{\cthr}{} & \cbox{white}{x_5} \\\hline}
    \end{array} 
  \end{displaymath}
    \caption{\label{fig:record}Skew diagrams for a saturated chain in $\kup_4$ from $\beta=(1,0,3,1,0,4)$ to $\gamma=(2,0,4,5,0,3)$ with $\alpha=(4,0,5,3,0,2) \lswap \gamma$.}
\end{figure}

For example, the chain in Fig.~\ref{fig:chain} corresponds to the skew diagrams in Fig.~\ref{fig:record}.

To insert a tableau into a diagram, we convert the tableau into a two-line array, where we insert the bottom line and use the top line to record landing cells. 

\begin{definition}
  For $T\in\KD(\beta)$ and $\left(\begin{array}{ccc}q_1 & \cdots & q_m \\ r_1 & \cdots & r_m \end{array}\right)$ a two-line array with each $r_i\le k$, the \newword{insertion diagram} is $U=(T \xleftarrow{\beta,k} r_1,\ldots,r_m)$ and the \newword{recording tableau} has shape $\thread(U)/\beta$ where the cell added when $r_i$ is inserted has entry $q_i$.
  \label{def:RSK-record}
\end{definition}

\begin{definition}
  For $\alpha\lswap\gamma$ and $\beta\kup_{k}\gamma$, an \newword{atomic tableau of skew shape $\alpha/\beta$} is a filling of the skew diagram $\alpha/\beta$ with positive integers such that
  \begin{enumerate}
  \item entries weakly decrease left to right within rows;
  \item entries within a column are distinct;
  \item if $t>r$ appear in the same column with $t$ above $r$, then there is an $s>r$ in the same row as and immediately to the right of $t$.
  \end{enumerate}
  Here we ignore cells above row $k$, labeled cells in column $1$ have entry $\infty$, and labeled cells weakly below $k$ have the same entry as their row neighbor to the left. 
  \label{def:AKT}
\end{definition}

The \newword{weight} has $\wt(T)_i$ equal to the number of skew cells with entry $i$. A tableau is \newword{lattice} if the weight of the columns weakly to the right of column $c$ is a partition. 

\begin{figure}[ht]
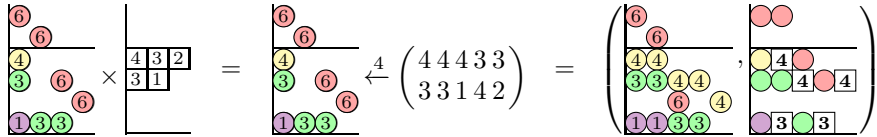

  \begin{displaymath}
    \arraycolsep=1pt
    \begin{array}{ccc}
      \cirtab{
        \cball{\csix}{6} \\ & \cball{\csix}{6} \\\hline
        \cball{\cfou}{4} \\
        \cball{\cthr}{3} & & \cball{\csix}{6} \\
        & & & \cball{\csix}{6} \\
        \cball{\cone}{1} & \cball{\cthr}{3} & \cball{\cthr}{3} &  \\\hline}
      & \raisebox{-2.7\cellsize}{$\times$} &
      \vline\tableau{ \\ \\ 4 & 3 & 2 \\ 3 & 1 \\ & \\ & \\\hline}
    \end{array}
    \hspace{\cellsize} = \hspace{\cellsize}
    \begin{array}{ccc}
      \cirtab{
        \cball{\csix}{6} \\ & \cball{\csix}{6} \\\hline
        \cball{\cfou}{4} \\
        \cball{\cthr}{3} & & \cball{\csix}{6} \\
        & & & \cball{\csix}{6} \\
        \cball{\cone}{1} & \cball{\cthr}{3} & \cball{\cthr}{3} &  \\\hline}
      & \raisebox{-2.7\cellsize}{$\xleftarrow{4}$} &
      \raisebox{-2.7\cellsize}{$\left(\begin{array}{ccccc}
        4 & 4 & 4 & 3 & 3 \\
        3 & 3 & 1 & 4 & 2 
      \end{array}\right)$}
    \end{array}
    \hspace{0.8\cellsize} = \hspace{0.8\cellsize}
    \left(
    \begin{array}{ccc}
      \vline\nulltab{
        \cball{\csix}{6} \\ & \cball{\csix}{6} \\\hline
        \cball{\cfou}{4} & \cball{\cfou}{4} & \\
        \cball{\cthr}{3} & \cball{\cthr}{3} & \cball{\cfou}{4} & \cball{\cfou}{4} \\
        & & \cball{\csix}{6} & & \cball{\cfou}{4} \\
        \cball{\cone}{1} & \cball{\cone}{1} & \cball{\cthr}{3} & \cball{\cthr}{3}  \\\hline}
      & \raisebox{-1.7\cellsize}{$,$} &
      \vline\nulltab{
        \cball{\csix}{} & \cball{\csix}{}   \\ \\\hline
        \cball{\cfou}{} & \cbox{white}{4} & \cball{\csix}{} \\ 
        \cball{\cthr}{} & \cball{\cthr}{} & \cbox{white}{4} & \cball{\csix}{} & \cbox{white}{4} \\
        & \\
        \cball{\cone}{} & \cbox{white}{3} & \cball{\cthr}{} & \cbox{white}{3} \\\hline}
    \end{array} \right)
  \end{displaymath}
  \caption{\label{fig:recording}The insertion diagram and atomic recording tableau for $T\in\AKD(10413)\subseteq\PKD(10314)\subseteq\KD(01314)$.}
\end{figure}

For $\lambda$ a partition, let $R_{\lambda}\in\SSRT(\lambda)$ be the tableau with every entry in row $r$ equal to $r$. Given $R\in\SSRT(\lambda)$, there is a unique two-line array with top row $(k^{\lambda_k},\ldots,1^{\lambda_1})$, where $k=\ell(\lambda)$, mapping to the pair $(R,R_{\lambda})$ under RSK (see \cite[\S4]{FYT}). The \newword{product} of a rectified diagram and a semistandard reverse tableau is the insertion of this two line array into the diagram; see Fig.~\ref{fig:recording}.

\begin{theorem}
  The recording tableau of the $k$-restricted insertion of a tableau $R\in\SSRT(\lambda)$ into a labeled diagram $T\in\KD(\beta)$ is a lattice atomic tableau of skew shape $\alpha/\beta$ for some $\alpha\lswap\gamma$ where $\beta\kup_{k}\gamma$ and weight $\lambda$.
  \label{thm:record}
\end{theorem}

\begin{proof}
  By Theorem~\ref{thm:iterate}, the map sending the pair $(T,R)\in\KD(\beta) \times \SSRT(\lambda)$ defined by iterated insertion of the word for $R$ into $T$ restricted by $k=\ell(\lambda)$ is well-defined and results in a diagram in $\AKD(\alpha)$ for some $\alpha\lswap\gamma$ where $\beta\kup_{k}\gamma$. Thus the recording tableau $Q$ is well-defined, and by definition, $\wt(Q)=\lambda$.

  Since the recorded entries weakly decrease and cells are added at the end of their rows, new entries within a row of $Q$ weakly decrease. Since cells are dropped only from above row $k$, the dropped cells do not change this property, proving condition (1) of Definition~\ref{def:AKT}. From properties of RSK (see \cite{FYT}), if two successive inserted letters, say $r,s$, have the same recording value, say $q$, then $r \ge s$. Thus by Lemma~\ref{lem:rowbump}\ref{i:dec}, the landing cells for the insertions lie in different columns, proving condition (2). Furthermore, by Definition~\ref{def:kadd}, if $r<t$ occur in the same column with $r$ below $t$ (which could be a cell with label at most $k$), then by Lemma~\ref{lem:rowbump}\ref{i:inc}, the cell $y$ corresponding to $r$ must have landed above the cell $x$ corresponding to $t$, 

  the row of $t$ cannot end in this column, else $r$ was not a $k$-extendable row. Therefore there is an entry to the right of $t$, say $s$, which was recorded prior to $r$ (or is labeled). Therefore $s>r$, and so condition (3) holds as well.

  For the lattice property, we need only show the column-restricted weights have at least as many $q$'s as they have $q-1$'s for all $q$. Suppose $r_1 \ge \cdots \ge r_m$ are recorded by $q$ and $s_1 \ge \cdots \ge s_n$ are recorded by $q-1$. Since RSK insertion of the two-line array records all $r_i$'s in row $q$ and all $s_i$'s in row $q-1$, we have $m \ge n$ and indices $i_1 < \ldots < i_n$ such that $r_{i_j} < s_{j}$. By Lemma~\ref{lem:rowbump}\ref{i:dec}, the bumping paths for the $r_i$'s are disjoint, and by Lemma~\ref{lem:rowbump}\ref{i:inc}, the landing cell for $s_j$ will lie weakly left of the landing cell for $r_{i_j}$. In particular, to the right of any column $c$, there are never more entries $q-1$ than entries $q$, so $Q$ is lattice.
\end{proof}

%
\section{Product rules for Demazure characters}
%
\label{sec:expansions}

In this section, we give nonnegative expansions into Demazure atoms and Schubert characters and lift these to \emph{signed} expansions into Demazure characters.

\subsection{Expansion into Demazure atoms}
\label{sec:2atoms}

All that remains is to account for the multiplicities using recording tableaux.

\begin{theorem}
  Restricted insertion gives a weight-preserving bijection
  \begin{equation}
    \KD(\beta) \times \SSRT(\lambda) \xrightarrow{\sim} \bigsqcup_{\substack{\beta\kup_{k}\gamma \\ \alpha\lswap\gamma}} \left( \AKD(\alpha) \times \LAT(\alpha/\beta,\lambda) \right),
    \label{e:atom-bijection}
  \end{equation}
  where $\LAT(\alpha/\beta,\lambda)$ is the set of lattice atomic tableaux of skew shape $\alpha/\beta$ and weight $\lambda$. In particular, we have a nonnegative expansion
  \begin{equation}
    \key_{\beta} s_{\lambda} = \sum_{\alpha} a_{\beta,\lambda}^{\alpha} \atom_{\alpha},
    \label{e:atom}
  \end{equation}
  where $a_{\beta,\lambda}^{\alpha}$ is the cardinality of $\LAT(\alpha/\beta,\lambda)$ and so is nonnegative.
  \label{thm:atom}
\end{theorem}

\begin{proof}
  By Theorem~\ref{thm:iterate}, the map sending the pair $(T,R)\in\KD(\beta) \times \SSRT(\lambda)$ defined by iterated insertion of the word for $R$ into $T$ restricted by $k=\ell(\lambda)$ is well-defined and results in a diagram in $\KD(\gamma)$ for some $\beta\kup_{k}\gamma$. By Lemma~\ref{lem:lswap}, $\AKD(\alpha)\subseteq\KD(\gamma)$ for all $\alpha\lswap\gamma$. By Theorem~\ref{thm:record}, the recording tableau is a lattice atomic tableau of skew shape $\alpha/\beta$ and weight $\lambda$. Moreover, the union for the image since threads are well-defined; that is, $\AKD(\alpha)\cap\AKD(\alpha') = \varnothing$ if $\alpha\neq\alpha'$.

  Given a pair $(U,Q)\in\AKD(\alpha)\times\LAT(\alpha/\beta,\lambda)$, by Lemma~\ref{lem:atomic} and \eqref{e:KD2AKD}, $U$ lies in the image of the map, though we must reverse with respect to $Q$. Let $Q' = Q\setminus\{x\}$ for $x$ the rightmost instance of the smallest entry of $Q$. Then $Q'$ is a lattice atomic tableau condition since $x$ is the smallest entry, necessarily taken from the end of a row. Let $c$ be the column of $x$ and $r$ the label in $\Label_{\gamma}(Q)$. Then there exists $\beta\kup_{k}\beta'\kupdot_{k}\gamma$ such that $\beta'\kplus(c,r)=\gamma$. By Lemma~\ref{lem:atomic}, there exists $T'\in\KD(\beta')$ and $r'\le k$ such that $U = (T'\xrightarrow{\beta',k}r')$. Moreover, by Theorem~\ref{thm:inverse}, both $T'$ and $r'$ are unique. By induction, the pair $(T',Q')$ uniquely corresponds to $(T,R)\in\KD(\beta)\times\SSRT(\lambda-\e_q)$, where $q$ is the entry of $x$. The theorem follows.
\end{proof}

\subsection{Expansion into Schubert characters}
\label{sec:2schubs}

The following construction allows us to lift the expansion into Demazure atoms in \eqref{e:atom} to Schubert characters.

\begin{definition}
  For $\beta\kup_{k}\gamma$, a \newword{key tableau of skew shape $\gamma/\beta$} is a filling of the skew diagram with positive integers such that
  \begin{enumerate}
  \item entries weakly decrease left to right within rows;
  \item entries within a column are distinct;
  \item if $t \ge r$ appear in adjacent columns with $r<\infty$ strictly below and right of $t$, then there is an entry $s>r$ immediately right of $t$.
  \end{enumerate}
  Here we ignore cells above row $k$, labeled cells in column $1$ have entry $\infty$, and labeled cells weakly below $k$ have the same entry as their neighbor to the left. 
  \label{def:YKT}
\end{definition}

Note for $\beta\kup_{k}\gamma$, an atomic tableau of shape $\gamma/\beta$ is usually not a key tableau. However, there is a bijection between atomic tableau and key tableau.


\begin{lemma}
  For $\alpha\lswap\gamma$ where $\beta\kup_{k}\gamma$, the map $\rkey_{\alpha}^{\gamma}:\LAT(\alpha/\beta)\rightarrow\LKT(\gamma/\beta)$ that assigns recorded entries to the added cells is a bijection.
  \label{lem:keytab}
\end{lemma}

\begin{proof}
  Consider $T\in\LAT(\alpha/\beta)$ and let $U=\rkey_{\alpha}^{\gamma}(T)$. Since cells are appended to the ends of rows, the rows of $U$ weakly decrease left to right. Since $\rkey_{\alpha}^{\gamma}$ preserves column sets, columns of $U$ have distinct entries and $U$ is lattice. In order for $(c,r)$ to be $k$-addable, for any row $r < r' \le k$, either there is no cell nor entry yet recorded in row $r'$, column $c-1$ or there is a cell or prior entry in row $r'$, column $c$ which necessarily has larger entry. Thus $U$ is a lattice key tableau.

  Conversely, let $U\in\LKT(\gamma/\beta)$. Let $(c,r)$ be the cell of the rightmost instance of the smallest entry, which is well-defined by condition (2). Then $r \le k$ and, by condition (1), $\gamma_r=c$, and by condition (3), $(c,r)$ is a $k$-addable cell for $\gamma-\e_r$, and so there exists $\beta\kup_{k}\beta'\kupdot_{k}\gamma$ such that $\beta'\kplus(c,r)$. Since $(c,r)$ has the smallest entry and lies at the end of its row, $U\setminus\{(c,r)\}$ is a lattice key tableau. Given any $\alpha\lswap\gamma$, by Theorem~\ref{thm:atom}, there exists a unique lattice atomic tableau $T$ of skew shape $\alpha/\beta'$ with the unique skew cell in column $c$. Let $\alpha'$ denote $\alpha$ with this cell removed. Then $\alpha'\lswap\beta'$, and by induction we are done.
\end{proof}

We have the following description of intersections of sets of Kohnert diagrams.

\begin{lemma}
  Given $\beta\kup_{k}\gamma^{(1)},\ldots,\gamma^{(m)}$ each with the same multiset of added columns taken in the same order, there exists $\gamma\lswap\gamma^{(1)},\ldots,\gamma^{(m)}$ such that
  \begin{equation}
    \KD(\gamma^{(1)}) \cap \cdots \cap \KD(\gamma^{(m)}) = \KD(\gamma).
  \end{equation}
  \label{lem:kap}
\end{lemma}

\begin{proof}
  By \cite[Lem~3.3.4]{AQ}, for $(c,r_1),\ldots,(c,r_m)$ $k$-addable positions for $\beta$, 
  \begin{equation}
    \bigcap_{i=1}^{m} \KD(\beta +_{k} (c,r_i)) = \KD(t_{r_1,r_2} \cdots t_{r_{m-1},r_m} (\beta +_{k} (c,r_m))),
    \label{e:kap}
  \end{equation}
  where indices are taken so that $r_1<\cdots<r_m$. Let $\beta'=t_{r_1,r_2} \cdots t_{r_{m-1},r_m} (\beta +_{k} (c,r_m))$. If $(c',r')$ is $k$-addable for all $\beta +_{k} (c,r_i)$, then it is $k$-addable for $\beta'$ as well. The result follows by induction on the length of the chain from $\beta$ to any $\gamma^{(i)}$. 
\end{proof}

We now prove Conjecture~\ref{conj:schubchar} for a Demazure character and Schur polynomial.

\begin{theorem}
  Given a sequence $\mathbf{c}=(c_1,\ldots,c_n)$ of successively $k$-addable columns for $\beta$, let $\Gamma_{\mathbf{c}}$ be the union of all compositions $\alpha\preceq\gamma$ for some $\gamma$ obtainable by successive $k$-addition of cells in these columns. Then $\Gamma_{\mathbf{c}}$ is a lower order ideal in Bruhat order, and we have a nonnegative expansion into Schubert characters as
  \begin{equation}
    \key_{\beta} s_{\lambda} = \sum_{\mathbf{c}} a_{\beta,\lambda}^{\Gamma_{\mathbf{c}}} \key_{\Gamma_{\mathbf{c}}},
    \label{e:schubcharexp}
  \end{equation}
  where $a_{\beta,\lambda}^{\Gamma_{\mathbf{c}}}$ is the cardinality of $\LAT(\alpha/\beta,\lambda)$ for any $\alpha\in\Gamma_{\mathbf{c}}$. 
  \label{thm:schubchar}
\end{theorem}

\begin{proof}
  It follows from the definition that $\alpha\in\Gamma_{\mathbf{c}}$ whenever $\alpha\preceq\gamma\in\Gamma_{\mathbf{c}}$, so $\Gamma_{\mathbf{c}}$ is a lower order ideal. By Lemma~\ref{lem:kap}, for any $\gamma,\gamma'\in\Gamma_{\mathbf{c}}$, there exists $\alpha\preceq\gamma,\gamma'$. Thus by Lemma~\ref{lem:keytab}, we have $\#\LAT(\gamma/\beta,\lambda) = \#\LAT(\alpha/\beta,\lambda) = \#\LAT(\gamma'/\beta,\lambda)$ for all $\lambda$. In particular, the cardinality of $\LAT(\alpha/\beta,\lambda)$ is the same for all $\alpha\in\Gamma_{\mathbf{c}}$. The result now follows from Theorem~\ref{thm:atom}.
\end{proof}

\subsection{Expansion into Demazure characters}
\label{sec:2keys}

Finally, we lift the expansion in \eqref{e:schubcharexp} to Demazure characters, though this introduces signs in many cases.

\begin{theorem}
  We have the cancellation-free, signed expansion
  \begin{equation}
    \key_{\beta} s_{\lambda} = \sum_{\beta\kup_{k}\gamma^{(1)},\ldots,\gamma^{(m)} \ \text{distinct}} (-1)^{m-1} c_{\beta,\lambda}^{\gamma^{(1)}} \key_{\gamma},
    \label{e:keyexp}
  \end{equation}
  where $\gamma$ is maximal such that $\gamma\lswap\gamma^{(i)}$ and $c_{\beta,\lambda}^{\gamma}$ is the cardinality of $\LKT(\gamma/\beta,\lambda)$, the set of lattice key tableaux of skew shape $\gamma/\beta$ and weight $\lambda$. 
  \label{thm:keyexp}
\end{theorem}

\begin{proof}
  By \eqref{e:KD2AKD} and Lemma~\ref{lem:keytab}, for $\beta\kup_{k}\gamma$ we have a bijection
  \[
  \bigsqcup_{\alpha\lswap\gamma} \AKD(\alpha)\times\LAT(\alpha/\beta,\lambda)
  \xrightarrow{\mathrm{id}\times\rkey_{\alpha}^{\gamma}}
  \KD(\gamma)\times\LKT(\gamma/\beta,\lambda).
  \]
  Expanding as in \eqref{e:atom-bijection} gives a union that is no longer disjoint. Thus \eqref{e:keyexp} follows from Lemma~\ref{lem:kap} by inclusion--exclusion.
\end{proof}

We can identify many instances in which the union in \eqref{e:atom-bijection} remains disjoint under the bijection in Lemma~\ref{lem:keytab}, and so the Demazure character expansion in \eqref{e:keyexp} is \emph{nonnegative}. 

\begin{definition}
  A composition $\beta$ is \newword{$k$-positive} if whenever $\beta_{r}\le \beta_{t}$ for some $r\le k < t$, then $\beta_s \ge \beta_r$ for all $r < s \le k$.
  \label{def:pos}
\end{definition}

\begin{theorem}
  For any $k$-positive $\beta$ (in particular, for $k \ge \ell(\beta)$), we have
  \begin{equation}
    \key_{\beta} s_{\lambda} = \sum_{\gamma} c_{\beta,\lambda}^{\gamma} \key_{\gamma},
    \label{e:keypos}
  \end{equation}
  where $c_{\beta,\lambda}^{\gamma}$ is the cardinality of $\LKT(\gamma/\beta,\lambda)$ and so is nonnegative.
  \label{thm:keypos}
\end{theorem}

\begin{proof}
  A composition $\beta$ has at most one $k$-addable position for each column if and only if whenever $\beta_{r}\le \beta_{t}$ for some $r\le k < t$, then either $\beta_s = \beta_r$ for some $r < s \le k$ or $\beta_s > \beta_r$ for all $r < s \le k$. Thus if $\beta$ is $k$-positive, then it has at most one $k$-addable position. Moreover, in this case, $\beta\kplus(c,r)$ is also $k$-positive for any $k$-addable position $(c,r)$. To see this, let $\gamma=\beta\kplus(c,r)$ and suppose, for contradiction, there exist $q < s \le k < t$ with $\gamma_s < \gamma_q \le \gamma_t$. Since $k$-addition lengthens row $r$ and shortens rows above $k$, leaving all others the same, in order for $\beta$ to be $k$-positive, we must have $q=r$ and $\beta_r < \beta_s < c$, where $t'$ is the maximum row from which cells were dropped. However, this contradicts that $(c,r)$ is $k$-addable. Thus we may $k$-add cells to $\beta$ uniquely based on the added column, making intersections in \eqref{e:keyexp} empty and the expansion nonnegative.
\end{proof}

In particular, since Kohnert diagrams naturally index a basis for the corresponding Demazure module, Theorem~\ref{thm:keypos} gives a large class of Demazure modules for which the tensor products admit excellent filtrations.

Many of these are not covered by recent work of Kouno \cite{Kou20} using Demazure crystals, since the crystals do not decompose into Demazure pieces despite the nonnegativity of the characters. That is, while RSK commutes with tensor product on crystals, the insertion algorithm constructed herein does not. Thus it would be interesting to explore if there is another tensor structure on crystals which does commute with this insertion.

%
%

\bibliographystyle{plain} 
\bibliography{atomic}

\begin{thebibliography}{10}

\bibitem{A-pnas}
Sami Assaf.
\newblock A conjectured {L}ittlewood--{R}ichardson rule for {G}rassmannian
  {S}chubert varieties.
\newblock {\em Proc. Natl. Acad. Sci.}
\newblock to appear.

\bibitem{AB}
Sami Assaf and Nantel Bergeron.
\newblock An insertion algorithm for multiplying {S}chubert polynomials by
  {S}chur polynomials.
\newblock 2023.

\bibitem{AQ}
Sami Assaf and Danjoseph Quijada.
\newblock Monk's rule for {D}emazure characters of the general linear group.
\newblock arXiv:1908.08502, 2019.

\bibitem{AS18}
Sami Assaf and Dominic Searles.
\newblock Kohnert tableaux and a lifting of quasi-{S}chur functions.
\newblock {\em J. Combin. Theory Ser. A}, 156:85--118, 2018.

\bibitem{Ass22-KC}
Sami~H. Assaf.
\newblock Demazure crystals for {K}ohnert polynomials.
\newblock {\em Trans. Amer. Math. Soc.}, 375(3):2147--2186, 2022.

\bibitem{Dem74a}
Michel Demazure.
\newblock D\'esingularisation des vari\'et\'es de {S}chubert
  g\'en\'eralis\'ees.
\newblock {\em Ann. Sci. \'Ecole Norm. Sup. (4)}, 7:53--88, 1974.
\newblock Collection of articles dedicated to Henri Cartan on the occasion of
  his 70th birthday, I.

\bibitem{Dem74}
Michel Demazure.
\newblock Une nouvelle formule des caract\`eres.
\newblock {\em Bull. Sci. Math. (2)}, 98(3):163--172, 1974.

\bibitem{FYT}
William Fulton.
\newblock {\em Young tableaux}, volume~35 of {\em London Mathematical Society
  Student Texts}.
\newblock Cambridge University Press, Cambridge, 1997.
\newblock With applications to representation theory and geometry.

\bibitem{HLMvW11}
J.~Haglund, K.~Luoto, S.~Mason, and S.~van Willigenburg.
\newblock Refinements of the {L}ittlewood-{R}ichardson rule.
\newblock {\em Trans. Amer. Math. Soc.}, 363(3):1665--1686, 2011.

\bibitem{Jos85}
A.~Joseph.
\newblock On the {D}emazure character formula.
\newblock {\em Ann. Sci. \'{E}cole Norm. Sup. (4)}, 18(3):389--419, 1985.

\bibitem{Knu70}
Donald~E. Knuth.
\newblock Permutations, matrices, and generalized {Y}oung tableaux.
\newblock {\em Pacific J. Math.}, 34:709--727, 1970.

\bibitem{Koh91}
Axel Kohnert.
\newblock Weintrauben, {P}olynome, {T}ableaux.
\newblock {\em Bayreuth. Math. Schr.}, (38):1--97, 1991.
\newblock Dissertation, Universit{\"a}t Bayreuth, Bayreuth, 1990.

\bibitem{Kou20}
Takafumi Kouno.
\newblock Decomposition of tensor products of {D}emazure crystals.
\newblock {\em J. Algebra}, 546:641--678, 2020.

\bibitem{LS90}
Alain Lascoux and Marcel-Paul Sch{\"u}tzenberger.
\newblock Keys \& standard bases.
\newblock In {\em Invariant theory and tableaux ({M}inneapolis, {MN}, 1988)},
  volume~19 of {\em IMA Vol. Math. Appl.}, pages 125--144. Springer, New York,
  1990.

\bibitem{Mas09}
Sarah Mason.
\newblock An explicit construction of type {A} {D}emazure atoms.
\newblock {\em J. Algebraic Combin.}, 29(3):295--313, 2009.

\bibitem{Mat89}
Olivier Mathieu.
\newblock Filtrations of {$B$}-modules.
\newblock {\em Duke Math. J.}, 59(2):421--442, 1989.

\bibitem{Pol89}
Patrick Polo.
\newblock Vari\'{e}t\'{e}s de {S}chubert et excellentes filtrations.
\newblock {\em Ast\'{e}risque}, (173-174):10--11, 281--311, 1989.
\newblock Orbites unipotentes et repr\'{e}sentations, III.

\bibitem{Pun16}
Anna~Ying Pun.
\newblock {\em On decomposition of the product of {D}emazure atoms and
  {D}emazure characters}.
\newblock ProQuest LLC, Ann Arbor, MI, 2016.
\newblock Thesis (Ph.D.)--University of Pennsylvania.

\bibitem{Rob38}
G.~de~B. Robinson.
\newblock On the {R}epresentations of the {S}ymmetric {G}roup.
\newblock {\em Amer. J. Math.}, 60(3):745--760, 1938.

\bibitem{Sch61}
C.~Schensted.
\newblock Longest increasing and decreasing subsequences.
\newblock {\em Canad. J. Math.}, 13:179--191, 1961.

\bibitem{Sea20}
Dominic Searles.
\newblock Polynomial bases: positivity and {S}chur multiplication.
\newblock {\em Trans. Amer. Math. Soc.}, 373(2):819--847, 2020.

\bibitem{vdK89}
Wilberd van~der Kallen.
\newblock Longest weight vectors and excellent filtrations.
\newblock {\em Math. Z.}, 201(1):19--31, 1989.

\end{thebibliography}

\appendix

\section{Example of $4$-insertion with bumping paths}

Fix $k=4$ and $\beta = (1,0,3,1,0,2)$. We will $4$-insert into rows $3,3,1,4,2$ of a Kohnert diagram for $\beta$. This corresponds to multiplying by the following insertion tableau under RSK since the recording tableau is lattice.
\[ \begin{array}{ccc}
  \left(\begin{array}{ccccc}
    4 & 4 & 4 & 3 & 3 \\
    3 & 3 & 1 & 4 & 2 
  \end{array}\right) &
  \xrightarrow{RSK} &
  \tableau{4 & 3 & 2 \\ 3 & 1} \times \tableau{4 & 4 & 4 \\ 3 & 3}
\end{array} \]
We indicate the bumping path and the atomic and key recording tableaux.

  \begin{displaymath}
    \arraycolsep=\cellsize
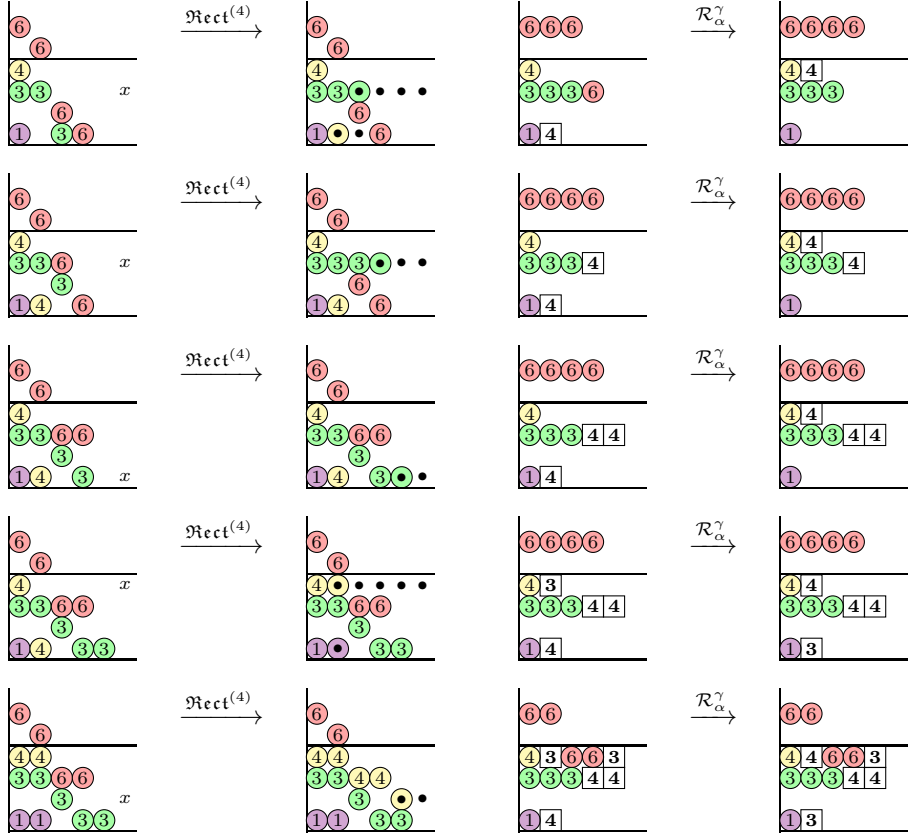
    \begin{array}{ccccccc}
      \vline\nulltab{
        \cball{\csix}{6} \\
        & \cball{\csix}{6} \\\hline
        \cball{\cfou}{4} & \\
        \cball{\cthr}{3} & \cball{\cthr}{3} & & & & x \\
        & & \cball{\csix}{6} \\
        \cball{\cone}{1} & & \cball{\cthr}{3} & \cball{\csix}{6}\\\hline }
      & \xrightarrow{\Rect^{(4)}} &
      \vline\nulltab{
        \cball{\csix}{6} \\
        & \cball{\csix}{6} \\\hline
        \cball{\cfou}{4} & \\
        \cball{\cthr}{3} & \cball{\cthr}{3} & \cball{\cthr}{\bullet} & \bullet & \bullet & \bullet \\
        & & \cball{\csix}{6} \\
        \cball{\cone}{1} & \cball{\cfou}{\bullet} & \bullet & \cball{\csix}{6}\\\hline }
      & &
      \vline\nulltab{
        \cball{\csix}{6} & \cball{\csix}{6} & \cball{\csix}{6} \\
        \\\hline
        \cball{\cfou}{4} & \\
        \cball{\cthr}{3} & \cball{\cthr}{3} & \cball{\cthr}{3} & \cball{\csix}{6} & & \\
        & & \\
        \cball{\cone}{1} & \cbox{white}{4} & & \\\hline }
      & \xrightarrow{\rkey_{\alpha}^{\gamma}} &
      \vline\nulltab{
        \cball{\csix}{6} & \cball{\csix}{6} & \cball{\csix}{6} & \cball{\csix}{6} \\
        \\\hline
        \cball{\cfou}{4} & \cbox{white}{4} \\
        \cball{\cthr}{3} & \cball{\cthr}{3} & \cball{\cthr}{3} & & & \\
        & & \\
        \cball{\cone}{1} & & & \\\hline }      
    \end{array}
  \end{displaymath}

  \begin{displaymath}
    \arraycolsep=\cellsize
    \begin{array}{ccccccc}
      \vline\nulltab{
        \cball{\csix}{6} \\
        & \cball{\csix}{6} \\\hline
        \cball{\cfou}{4} & \\
        \cball{\cthr}{3} & \cball{\cthr}{3} & \cball{\csix}{6} &  &  & x \\
        & & \cball{\cthr}{3} \\
        \cball{\cone}{1} & \cball{\cfou}{4} & & \cball{\csix}{6}\\\hline }
      & \xrightarrow{\Rect^{(4)}} &
      \vline\nulltab{
        \cball{\csix}{6} \\
        & \cball{\csix}{6} \\\hline
        \cball{\cfou}{4} & \\
        \cball{\cthr}{3} & \cball{\cthr}{3} & \cball{\cthr}{3} & \cball{\cthr}{\bullet} & \bullet & \bullet \\
        & & \cball{\csix}{6} \\
        \cball{\cone}{1} & \cball{\cfou}{4} & & \cball{\csix}{6}\\\hline }
      & &
      \vline\nulltab{
        \cball{\csix}{6} & \cball{\csix}{6} & \cball{\csix}{6} & \cball{\csix}{6} \\
        \\\hline
        \cball{\cfou}{4} & \\
        \cball{\cthr}{3} & \cball{\cthr}{3} & \cball{\cthr}{3} & \cbox{white}{4} & & \\
        & & \\
        \cball{\cone}{1} & \cbox{white}{4} & & \\\hline }
      & \xrightarrow{\rkey_{\alpha}^{\gamma}} &
      \vline\nulltab{
        \cball{\csix}{6} & \cball{\csix}{6} & \cball{\csix}{6} & \cball{\csix}{6} \\
        \\\hline
        \cball{\cfou}{4} & \cbox{white}{4} \\
        \cball{\cthr}{3} & \cball{\cthr}{3} & \cball{\cthr}{3} & \cbox{white}{4} & & \\
        & & \\
        \cball{\cone}{1} & & & \\\hline }      
    \end{array}
  \end{displaymath}

    \begin{displaymath}
    \arraycolsep=\cellsize
    \begin{array}{ccccccc}
      \vline\nulltab{
        \cball{\csix}{6} \\
        & \cball{\csix}{6} \\\hline
        \cball{\cfou}{4} & \\
        \cball{\cthr}{3} & \cball{\cthr}{3} & \cball{\csix}{6} & \cball{\csix}{6} & &  \\
        & & \cball{\cthr}{3} \\
        \cball{\cone}{1} & \cball{\cfou}{4} & & \cball{\cthr}{3} & & x\\\hline }
      & \xrightarrow{\Rect^{(4)}} &
      \vline\nulltab{
        \cball{\csix}{6} \\
        & \cball{\csix}{6} \\\hline
        \cball{\cfou}{4} & \\
        \cball{\cthr}{3} & \cball{\cthr}{3} & \cball{\csix}{6} & \cball{\csix}{6} & &  \\
        & & \cball{\cthr}{3} \\
        \cball{\cone}{1} & \cball{\cfou}{4} & & \cball{\cthr}{3} & \cball{\cthr}{\bullet} & \bullet \\\hline }
      & &
      \vline\nulltab{
        \cball{\csix}{6} & \cball{\csix}{6} & \cball{\csix}{6} & \cball{\csix}{6} \\
        \\\hline
        \cball{\cfou}{4} & \\
        \cball{\cthr}{3} & \cball{\cthr}{3} & \cball{\cthr}{3} & \cbox{white}{4} & \cbox{white}{4} & \\
        & & \\
        \cball{\cone}{1} & \cbox{white}{4} & & \\\hline }
      & \xrightarrow{\rkey_{\alpha}^{\gamma}} &
      \vline\nulltab{
        \cball{\csix}{6} & \cball{\csix}{6} & \cball{\csix}{6} & \cball{\csix}{6} \\
        \\\hline
        \cball{\cfou}{4} & \cbox{white}{4} \\
        \cball{\cthr}{3} & \cball{\cthr}{3} & \cball{\cthr}{3} & \cbox{white}{4} & \cbox{white}{4} & \\
        & & \\
        \cball{\cone}{1} & & & \\\hline }      
    \end{array}
  \end{displaymath}

    \begin{displaymath}
    \arraycolsep=\cellsize
    \begin{array}{ccccccc}
      \vline\nulltab{
        \cball{\csix}{6} \\
        & \cball{\csix}{6} \\\hline
        \cball{\cfou}{4} & & & & & x \\
        \cball{\cthr}{3} & \cball{\cthr}{3} & \cball{\csix}{6} & \cball{\csix}{6} & &  \\
        & & \cball{\cthr}{3} \\
        \cball{\cone}{1} & \cball{\cfou}{4} & & \cball{\cthr}{3} & \cball{\cthr}{3} &  \\\hline }
      & \xrightarrow{\Rect^{(4)}} &
      \vline\nulltab{
        \cball{\csix}{6} \\
        & \cball{\csix}{6} \\\hline
        \cball{\cfou}{4} & \cball{\cfou}{\bullet} & \bullet & \bullet & \bullet & \bullet \\
        \cball{\cthr}{3} & \cball{\cthr}{3} & \cball{\csix}{6} & \cball{\csix}{6} & &  \\
        & & \cball{\cthr}{3} \\
        \cball{\cone}{1} & \cball{\cone}{\bullet} & & \cball{\cthr}{3} & \cball{\cthr}{3} &  \\\hline }
      & &
      \vline\nulltab{
        \cball{\csix}{6} & \cball{\csix}{6} & \cball{\csix}{6} & \cball{\csix}{6} \\
        \\\hline
        \cball{\cfou}{4} & \cbox{white}{3} \\
        \cball{\cthr}{3} & \cball{\cthr}{3} & \cball{\cthr}{3} & \cbox{white}{4} & \cbox{white}{4} & \\
        & & \\
        \cball{\cone}{1} & \cbox{white}{4} & & \\\hline }
      & \xrightarrow{\rkey_{\alpha}^{\gamma}} &
      \vline\nulltab{
        \cball{\csix}{6} & \cball{\csix}{6} & \cball{\csix}{6} & \cball{\csix}{6} \\
        \\\hline
        \cball{\cfou}{4} & \cbox{white}{4} \\
        \cball{\cthr}{3} & \cball{\cthr}{3} & \cball{\cthr}{3} & \cbox{white}{4} & \cbox{white}{4} & \\
        & & \\
        \cball{\cone}{1} & \cbox{white}{3} & & \\\hline }      
    \end{array}
    \end{displaymath}

    \begin{displaymath}
    \arraycolsep=\cellsize
    \begin{array}{ccccccc}
      \vline\nulltab{
        \cball{\csix}{6} \\
        & \cball{\csix}{6} \\\hline
        \cball{\cfou}{4} & \cball{\cfou}{4} \\
        \cball{\cthr}{3} & \cball{\cthr}{3} & \cball{\csix}{6} & \cball{\csix}{6} & &  \\
        & & \cball{\cthr}{3} & & & x \\
        \cball{\cone}{1} & \cball{\cone}{1} & & \cball{\cthr}{3} & \cball{\cthr}{3} & \\\hline }
      & \xrightarrow{\Rect^{(4)}} &
      \vline\nulltab{
        \cball{\csix}{6} \\
        & \cball{\csix}{6} \\\hline
        \cball{\cfou}{4} & \cball{\cfou}{4} \\
        \cball{\cthr}{3} & \cball{\cthr}{3} & \cball{\cfou}{4} & \cball{\cfou}{4} & &  \\
        & & \cball{\cthr}{3} & & \cball{\cfou}{\bullet} & \bullet \\
        \cball{\cone}{1} & \cball{\cone}{1} & & \cball{\cthr}{3} & \cball{\cthr}{3} & \\\hline }
      & &
      \vline\nulltab{
        \cball{\csix}{6} & \cball{\csix}{6}  \\
        \\\hline
        \cball{\cfou}{4} & \cbox{white}{3} & \cball{\csix}{6} & \cball{\csix}{6} & \cbox{white}{3} \\
        \cball{\cthr}{3} & \cball{\cthr}{3} & \cball{\cthr}{3} & \cbox{white}{4} & \cbox{white}{4} & \\
        & & \\
        \cball{\cone}{1} & \cbox{white}{4} & & \\\hline }
      & \xrightarrow{\rkey_{\alpha}^{\gamma}} &
      \vline\nulltab{
        \cball{\csix}{6} & \cball{\csix}{6}  \\
        \\\hline
        \cball{\cfou}{4} & \cbox{white}{4} & \cball{\csix}{6} & \cball{\csix}{6} & \cbox{white}{3} \\
        \cball{\cthr}{3} & \cball{\cthr}{3} & \cball{\cthr}{3} & \cbox{white}{4} & \cbox{white}{4} & \\
        & & \\
        \cball{\cone}{1} & \cbox{white}{3} & & \\\hline }      
    \end{array}
    \end{displaymath}

\section{Example of $4$-extraction with bumping paths}

Fix $k=4$ and the $4$-chain from $\beta = (1,0,3,1,0,4)$ to $\gamma = (2,0,4,5,0,3)$ from Fig.~\ref{fig:chain}. We will reverse the insertion in columns $4,2,5,3,2$ to recover the original Kohnert diagram for $\beta$ into which five rows were inserted.

As before, we give the bumping path and the atomic and key recording tableaux.

    \begin{displaymath}
    \arraycolsep=\cellsize
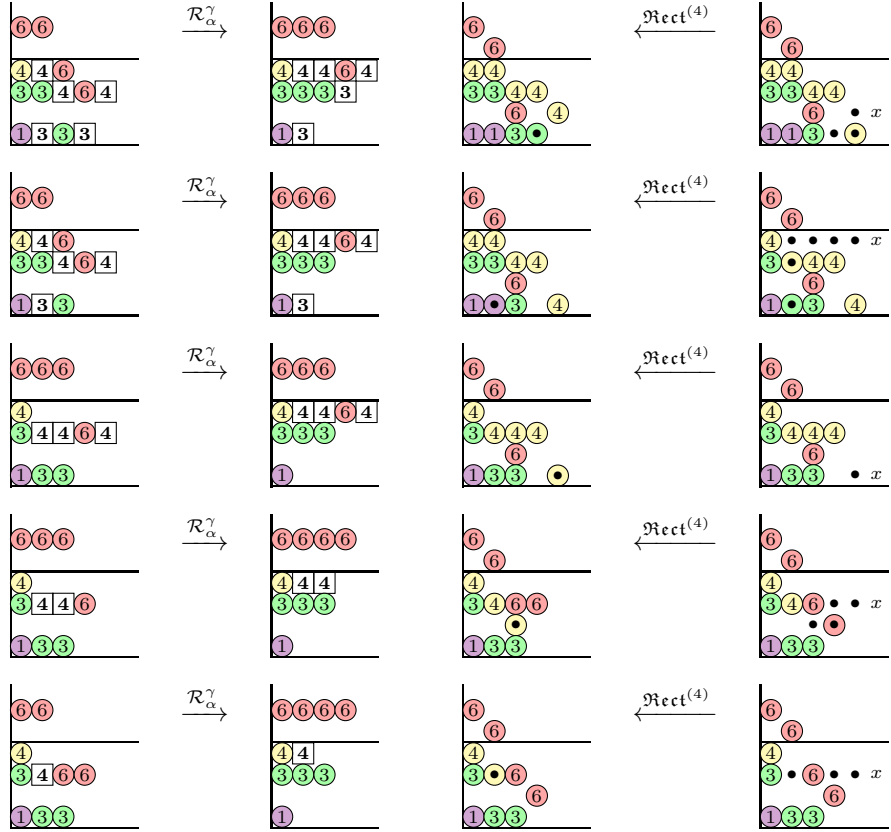
    \begin{array}{ccccccc}
      \vline\nulltab{
        \cball{\csix}{6} & \cball{\csix}{6}  \\
        \\\hline
        \cball{\cfou}{4} & \cbox{white}{4} & \cball{\csix}{6}  \\
        \cball{\cthr}{3} & \cball{\cthr}{3} & \cbox{white}{4} & \cball{\csix}{6} & \cbox{white}{4} & \\
        & & \\
        \cball{\cone}{1} & \cbox{white}{3} & \cball{\cthr}{3} & \cbox{white}{3} \\\hline }
      & \xrightarrow{\rkey_{\alpha}^{\gamma}} &
      \vline\nulltab{
        \cball{\csix}{6} & \cball{\csix}{6} & \cball{\csix}{6}  \\
        \\\hline
        \cball{\cfou}{4} & \cbox{white}{4} & \cbox{white}{4} & \cball{\csix}{6} & \cbox{white}{4} \\
        \cball{\cthr}{3} & \cball{\cthr}{3} & \cball{\cthr}{3} & \cbox{white}{3} & \\
        & & \\
        \cball{\cone}{1} & \cbox{white}{3} & & \\\hline }
      & &
      \vline\nulltab{
        \cball{\csix}{6} \\
        & \cball{\csix}{6} \\\hline
        \cball{\cfou}{4} & \cball{\cfou}{4} \\
        \cball{\cthr}{3} & \cball{\cthr}{3} & \cball{\cfou}{4} & \cball{\cfou}{4} & &  \\
        & & \cball{\csix}{6} & & \cball{\cfou}{4} &  \\
        \cball{\cone}{1} & \cball{\cone}{1} & \cball{\cthr}{3} & \cball{\cthr}{\bullet} & & \\\hline }
      & \xleftarrow{\Rect^{(4)}} &
      \vline\nulltab{
        \cball{\csix}{6} \\
        & \cball{\csix}{6} \\\hline
        \cball{\cfou}{4} & \cball{\cfou}{4} \\
        \cball{\cthr}{3} & \cball{\cthr}{3} & \cball{\cfou}{4} & \cball{\cfou}{4} & &  \\
        & & \cball{\csix}{6} & & \bullet & x \\
        \cball{\cone}{1} & \cball{\cone}{1} & \cball{\cthr}{3} & \bullet & \cball{\cfou}{\bullet} & \\\hline }
    \end{array}
    \end{displaymath}

    \begin{displaymath}
    \arraycolsep=\cellsize
    \begin{array}{ccccccc}
      \vline\nulltab{
        \cball{\csix}{6} & \cball{\csix}{6}  \\
        \\\hline
        \cball{\cfou}{4} & \cbox{white}{4} & \cball{\csix}{6}  \\
        \cball{\cthr}{3} & \cball{\cthr}{3} & \cbox{white}{4} & \cball{\csix}{6} & \cbox{white}{4} & \\
        & & \\
        \cball{\cone}{1} & \cbox{white}{3} & \cball{\cthr}{3}  \\\hline }
      & \xrightarrow{\rkey_{\alpha}^{\gamma}} &
      \vline\nulltab{
        \cball{\csix}{6} & \cball{\csix}{6} & \cball{\csix}{6}  \\
        \\\hline
        \cball{\cfou}{4} & \cbox{white}{4} & \cbox{white}{4} & \cball{\csix}{6} & \cbox{white}{4} \\
        \cball{\cthr}{3} & \cball{\cthr}{3} & \cball{\cthr}{3} & & \\
        & & \\
        \cball{\cone}{1} & \cbox{white}{3} & & \\\hline }
      & &
      \vline\nulltab{
        \cball{\csix}{6} \\
        & \cball{\csix}{6} \\\hline
        \cball{\cfou}{4} & \cball{\cfou}{4} \\
        \cball{\cthr}{3} & \cball{\cthr}{3} & \cball{\cfou}{4} & \cball{\cfou}{4} & &  \\
        & & \cball{\csix}{6} & &  &  \\
        \cball{\cone}{1} & \cball{\cone}{\bullet} & \cball{\cthr}{3} & & \cball{\cfou}{4} & \\\hline }
      & \xleftarrow{\Rect^{(4)}} &
      \vline\nulltab{
        \cball{\csix}{6} \\
        & \cball{\csix}{6} \\\hline
        \cball{\cfou}{4} & \bullet & \bullet & \bullet & \bullet & x \\
        \cball{\cthr}{3} & \cball{\cfou}{\bullet} & \cball{\cfou}{4} & \cball{\cfou}{4} & &  \\
        & & \cball{\csix}{6} & &  &  \\
        \cball{\cone}{1} & \cball{\cthr}{\bullet} & \cball{\cthr}{3} & & \cball{\cfou}{4} & \\\hline }
    \end{array}
    \end{displaymath}

    \begin{displaymath}
    \arraycolsep=\cellsize
    \begin{array}{ccccccc}
      \vline\nulltab{
        \cball{\csix}{6} & \cball{\csix}{6} & \cball{\csix}{6} \\
        \\\hline
        \cball{\cfou}{4} &  \\
        \cball{\cthr}{3} & \cbox{white}{4} & \cbox{white}{4} & \cball{\csix}{6} & \cbox{white}{4} & \\
        & & \\
        \cball{\cone}{1} & \cball{\cthr}{3} & \cball{\cthr}{3}  \\\hline }
      & \xrightarrow{\rkey_{\alpha}^{\gamma}} &
      \vline\nulltab{
        \cball{\csix}{6} & \cball{\csix}{6} & \cball{\csix}{6}  \\
        \\\hline
        \cball{\cfou}{4} & \cbox{white}{4} & \cbox{white}{4} & \cball{\csix}{6} & \cbox{white}{4} \\
        \cball{\cthr}{3} & \cball{\cthr}{3} & \cball{\cthr}{3} & & \\
        & & \\
        \cball{\cone}{1} & & & \\\hline }
      & &
      \vline\nulltab{
        \cball{\csix}{6} \\
        & \cball{\csix}{6} \\\hline
        \cball{\cfou}{4} &  \\
        \cball{\cthr}{3} & \cball{\cfou}{4} & \cball{\cfou}{4} & \cball{\cfou}{4} & &  \\
        & & \cball{\csix}{6} & &  &  \\
        \cball{\cone}{1} & \cball{\cthr}{3} & \cball{\cthr}{3} & & \cball{\cfou}{\bullet} & \\\hline }
      & \xleftarrow{\Rect^{(4)}} &
      \vline\nulltab{
        \cball{\csix}{6} \\
        & \cball{\csix}{6} \\\hline
        \cball{\cfou}{4} &  \\
        \cball{\cthr}{3} & \cball{\cfou}{4} & \cball{\cfou}{4} & \cball{\cfou}{4} & &  \\
        & & \cball{\csix}{6} & &  &  \\
        \cball{\cone}{1} & \cball{\cthr}{3} & \cball{\cthr}{3} & & \bullet & x \\\hline }
    \end{array}
    \end{displaymath}

    \begin{displaymath}
    \arraycolsep=\cellsize
    \begin{array}{ccccccc}
      \vline\nulltab{
        \cball{\csix}{6} & \cball{\csix}{6} & \cball{\csix}{6} \\
        \\\hline
        \cball{\cfou}{4} &  \\
        \cball{\cthr}{3} & \cbox{white}{4} & \cbox{white}{4} & \cball{\csix}{6} & & \\
        & & \\
        \cball{\cone}{1} & \cball{\cthr}{3} & \cball{\cthr}{3}  \\\hline }
      & \xrightarrow{\rkey_{\alpha}^{\gamma}} &
      \vline\nulltab{
        \cball{\csix}{6} & \cball{\csix}{6} & \cball{\csix}{6} & \cball{\csix}{6} \\
        \\\hline
        \cball{\cfou}{4} & \cbox{white}{4} & \cbox{white}{4} &  & \\
        \cball{\cthr}{3} & \cball{\cthr}{3} & \cball{\cthr}{3} & & \\
        & & \\
        \cball{\cone}{1} & & & \\\hline }
      & &
      \vline\nulltab{
        \cball{\csix}{6} \\
        & \cball{\csix}{6} \\\hline
        \cball{\cfou}{4} &  \\
        \cball{\cthr}{3} & \cball{\cfou}{4} & \cball{\csix}{6} & \cball{\csix}{6} & &  \\
        & & \cball{\cfou}{\bullet} & &  &  \\
        \cball{\cone}{1} & \cball{\cthr}{3} & \cball{\cthr}{3} & & & \\\hline }
      & \xleftarrow{\Rect^{(4)}} &
      \vline\nulltab{
        \cball{\csix}{6} \\
        & \cball{\csix}{6} \\\hline
        \cball{\cfou}{4} &  \\
        \cball{\cthr}{3} & \cball{\cfou}{4} & \cball{\csix}{6} & \bullet & \bullet & x \\
        & & \bullet & \cball{\csix}{\bullet} &  &  \\
        \cball{\cone}{1} & \cball{\cthr}{3} & \cball{\cthr}{3} & & & \\\hline }
    \end{array}
    \end{displaymath}

    \begin{displaymath}
    \arraycolsep=\cellsize
    \begin{array}{ccccccc}
      \vline\nulltab{
        \cball{\csix}{6} & \cball{\csix}{6} &  \\
        \\\hline
        \cball{\cfou}{4} &  \\
        \cball{\cthr}{3} & \cbox{white}{4} & \cball{\csix}{6} & \cball{\csix}{6} & & \\
        & & \\
        \cball{\cone}{1} & \cball{\cthr}{3} & \cball{\cthr}{3}  \\\hline }
      & \xrightarrow{\rkey_{\alpha}^{\gamma}} &
      \vline\nulltab{
        \cball{\csix}{6} & \cball{\csix}{6} & \cball{\csix}{6} & \cball{\csix}{6} \\
        \\\hline
        \cball{\cfou}{4} & \cbox{white}{4} & &  & \\
        \cball{\cthr}{3} & \cball{\cthr}{3} & \cball{\cthr}{3} & & \\
        & & \\
        \cball{\cone}{1} & & & \\\hline }
      & & 
      \vline\nulltab{
        \cball{\csix}{6} \\
        & \cball{\csix}{6} \\\hline
        \cball{\cfou}{4} &  \\
        \cball{\cthr}{3} & \cball{\cfou}{\bullet} & \cball{\csix}{6} & & &  \\
        & & & \cball{\csix}{6} & &  \\
        \cball{\cone}{1} & \cball{\cthr}{3} & \cball{\cthr}{3} & & & \\\hline }
      & \xleftarrow{\Rect^{(4)}} &
      \vline\nulltab{
        \cball{\csix}{6} \\
        & \cball{\csix}{6} \\\hline
        \cball{\cfou}{4} &  \\
        \cball{\cthr}{3} & \bullet & \cball{\csix}{6} & \bullet & \bullet & x \\
        & & & \cball{\csix}{6} & &  \\
        \cball{\cone}{1} & \cball{\cthr}{3} & \cball{\cthr}{3} & & & \\\hline }
    \end{array}
    \end{displaymath}

    Thus the extended rows were $2,4,1,3,3$. In particular, we have verified the insertion in Fig.~\ref{fig:recording} by excision.

\end{document}